\documentclass[12 pt]{amsart}
\usepackage{amsmath}
\usepackage{amsfonts}
\usepackage{amsthm}
\usepackage{amssymb}
\usepackage{amscd}
\usepackage[all]{xy}
\usepackage{enumerate}

\title{The moduli space of genus $4$ even spin curves is rational}

\author{Hiromichi Takagi and Francesco Zucconi}
\address{
Graduate School of Mathematical Sciences, 
the University of Tokyo,
Tokyo, 153-8914, Japan D.I.M.I. Universit\`a di Udine Via delle
scienze
208-206,
33100 Udine Italia}
\email{takagi@ms.u-tokyo.ac.jp
zucconi@dimi.uniud.it}

\theoremstyle{plain}
\newtheorem{thm}{Theorem}[section]

\newtheorem{prop}[thm]{Proposition}

\newtheorem{cor}[thm]{Corollary}
\newtheorem{lem}[thm]{Lemma}
\newtheorem{cla}[thm]{Claim}

\theoremstyle{definition}
\newtheorem{defn}[thm]{Definition}

\newtheorem{expl}[thm]{Example}

\newtheorem{nota}[thm]{Notation}

\newtheorem*{ackn}{Acknowledgment}

\theoremstyle{remark}
\newtheorem*{rem}{Remark}

\newcommand{\sD}{{\mathcal D}}
\newcommand{\sE}{{\mathcal E}}
\newcommand{\sF}{{\mathcal F}}
\newcommand{\sG}{{\mathcal G}}
\newcommand{\sH}{{\mathcal H}}

\newcommand{\sL}{{\mathcal L}}
\newcommand{\sN}{{\mathcal N}}
\newcommand{\sM}{{\mathcal M}}
\newcommand{\sO}{{\mathcal O}}

\newcommand{\sR}{{\mathcal R}}
\newcommand{\sS}{{\mathcal S}}
\newcommand{\sT}{{\mathcal T}}

\newcommand{\sV}{{\mathcal V}}
\newcommand{\sW}{{\mathcal W}}

\newcommand{\sY}{{\mathcal Y}}

\newcommand{\mC}{{\mathbb C}}

\newcommand{\mN}{{\mathbb N}}
\newcommand{\mP}{{\mathbb P}}
\newcommand{\mQ}{{\mathbb Q}}
\newcommand{\mZ}{{\mathbb Z}}

\newcommand{\Ima}{\mathrm{Im}\,}

\newcommand{\Bs}{\mathrm{Bs}\,}

\newcommand{\Aut}{\mathrm{Aut}\,}
\newcommand{\Hilb}{\mathrm{Hilb}}

\newcommand{\Pic}{\mathrm{Pic}\,}

\newcommand{\PGL}{\mathrm{PGL}}

\newcommand{\Supp}{\mathrm{Supp}\,}

\newcommand{\VSP}{\mathrm{VSP}\,}

\numberwithin{equation}{section}

\begin{document}


\maketitle
\markboth{Takagi and Zucconi}{Spin curves of genus $4$}

\begin{abstract}
    By the technique of $3$-fold Mori theory,
    we prove that the moduli space whose general point parameterizes a
    couple $(\sH,\theta)$ of a smooth curve $\sH$ of genus $4$ and
    a halfcanonical divisor $\theta$ such that 
    $h^{0}(\sH,\sO_{\sH}(\theta))=0$ is birational to $\mP^{9}$.
\end{abstract}



\section{Introduction}
Throughout the paper, we work over $\mC$, the complex number field.

A {\it{spin curve}} is a couple $(\Gamma,\theta)$,
where $\Gamma$ is a smooth projective curve of genus $g$
and $\theta$ is a {\em{theta characteristic}},
that is, an element $\theta\in {\rm{Pic}}\,\Gamma$ 
such that $2\theta$ is the class of the canonical sheaf $\omega_{\Gamma}$.
There are $2^{2g}$ different kinds of spin curve structures
for every smooth curve $\Gamma$ and they are partitioned into two classes according to the parity of $h^{0}(\Gamma,\theta)$.
A theta characteristic $\theta$ is said to be {\em{even}} or {\em{odd}}
if $h^{0}(\Gamma,\theta)$ is even or odd respectively.
Correspondingly we speak of
{\em{even or odd spin curves}}.

There exists the moduli space $\sS_{g}$ which parameterizes 
smooth spin curves $(\Gamma,\theta)$
and by the forgetful map $\sS_{g}\rightarrow \sM_{g}$,
where $\sM_{g}$ is the moduli space of curves of genus $g$,
we see that $\sS_{g}$ is a disjoint union of two irreducible
components $\sS^{+}_{g}$ and $\sS^{-}_{g}$ of relative degrees
$2^{g-1}(2^{g}+1)$ and $2^{g-1}(2^{g}-1)$ corresponding to even and
odd spin curves respectively (\cite{Mum1}, \cite{ACGH}).  
It was classically known that
$\sS^{+}_2$ is rational. 
The so called Scorza map gives
a birational isomorphisms between $\sS^{+}_3$ and
$\sM_3$ (\cite{DK}), hence $\sS^+_3$ is rational since
so is $\sM_3$ by \cite{Kat} (see also \cite{Boh}). 
Recently, Farkas \cite{Farkas} proved that 
a compactification $\overline{\sS}^+_g$ of ${\sS}^+_g$ 
is of general type for $g>8$, and  
the Kodaira dimension of $\overline{\sS}^+_g$ is negative 
for $g<8$, and non-negative for $g=8$.

In the previous papers \cite{TZ} and \cite{TZb},
we discovered a method to study trigonal even spin curves of any genus
by using biregular and birational geometries of the quintic del Pezzo $3$-fold
$B$.
The $3$-fold $B$ is, by definition, a smooth projective threefold
such that $-K_{B}=2H$, where $H$ is the ample generator
of $\Pic B$ and $H^3=5$. It is well known that the linear system 
$|H|$ embeds $B$ into $\mP^6$.

We explain our method specializing to the genus $4$ case,
which is under consideration in this paper. 
In this case, our main ingredient is  
the Hilbert scheme $\sH$ of 
general sextic normal rational curves on $B$.
We have shown that $\sH$ is irreducible
(see \cite[Proposition 2.5.2]{TZ}).
For a general sextic normal rational curve $C$ on $B$,
we have constructed a smooth curve $\sH_1$ of genus $4$
and a theta characteristic $\theta$ on it.
They come from the geometry of lines on $B$ intersecting $C$.
It is known that $\Aut B$ is isomorphic to
the automorphism group $\mathrm{PGL}_2$ of 
the complex projective line (see \cite{MU} and \cite{PV}).
From now on we set $G:=\PGL_2$.
The $G$-action on $B$ induces a $G$-action on $\sH$.
Thus we have a natural rational map
$\pi_{\sS^{+}_{4}}\colon  
\sH\dashrightarrow \sS^{+}_{4}$
which maps a general $C$ to $(\sH_1,\theta)$
and is constant on general $G$-orbits.
By taking suitable compactifications of $\sH$ and of
$\sS^{+}_{4}$,
a resolution of indeterminancy of $\pi_{\sS^{+}_{4}}$ and
the Stein factorization,
we have rational maps
$p_{\sS^{+}_{4}}\colon \sH\dashrightarrow
\widetilde{\sS}^{+}_{4}$
and $q_{\sS^{+}_{4}}\colon \widetilde{\sS}^{+}_{4}\dashrightarrow
{\sS}^{+}_{4}$ 
such that $\pi_{\sS^{+}_{4}}$ is given by $q_{\sS^{+}_{4}}\circ
p_{\sS^{+}_{4}}$, a general fiber of $p_{\sS^{+}_{4}}$ is connected and
$q_{\sS^{+}_{4}}$ is generically finite.
Then the $G$-orbit of a general point of $\sH$
is contained in a fiber of $p_{\sS^{+}_{4}}$.
In \cite[Theorem 4.0.2]{TZb}, we have proved that 
$\widetilde{\sS}^+_{4}$ is birational to $\sS^+_4$
or to its double cover, and 
birationally parameterizes $G$-orbits in $\sH$. 

Farkas's result mentioned above and 
the rationality of $\sM_4$ (\cite{SShB}) motivated us
to deepen our understanding of ${\sS}^{+}_{4}$. 
Then we obtain the following result in this paper:

\begin{thm}
$\sS^+_4$ is rational.
\end{thm}

Roughly speaking, the paper essentially consists of three parts;
in the section \ref{section ben noto}, we review the results 
as for the biregular geometries of $B$.
Especially, we review the biregular descriptions of $B$ and
the behaviour of lines on $B$ and sextic rational curves on $B$.
We also review the construction of the even spin curve $(\sH_1,\theta)$
of genus $4$ from a sextic rational curve on $B$.
In the section 3, 
we study some special birational selfmap $B\dashrightarrow B$,
which is one of our main ingredients to show 
the rationality of $\sS^+_4$.
Finally, in the section 4, we prove the rationality of $\sS^+_4$
as applications of the results in the section 2 and 3.

Now we explain an outline of the proof of the rationality of $\sS^+_4$.
Among other things in the section \ref{section ben noto}, 
we remind the readers 
that a general sextic normal smooth rational curve on $B$
has a natural marking, namely, 
its $6$ distinct bi-secant lines on $B$ (see Corollary \ref{orabisecanti}).
Therefore, noting the Hilbert scheme $\sH^B_1$ of lines on $B$ is $\mP^2$
(see the subsection \ref{quasinoto}),
we can define the morphism 
$\Theta\colon \widetilde{U}_0\to (\mP^2)^{6}/\mathfrak{S}_6$
mapping a sextic curve
to the unordered set of its $6$ bi-secant lines, 
where $\widetilde{U}_0$ is the open subset of $\sH$ consisting of 
sextic curves with exactly six different bi-secant lines, and
$(\mP^2)^{6}$
is the Cartesian product of six copies of $\mP^{2}\simeq \sH^{B}_{1}$ and 
$\mathfrak{S}_6$
is the permutation group on its factors. 
In Theorem \ref{pointsonly},
we show that $\Theta$ is birational.

Its proof requires the detailed study presented 
in the section 3 of the above mentioned
birational selfmap $B\dashrightarrow B$ centered 
along a smooth sextic rational curve
(Proposition \ref{prop:tworay}). 
There we use techniques of the $3$-fold explicit Mori theory,
especially, properties of smooth flops and the classification of
extremal contractions from smooth $3$-folds.
The selfmap $B\dashrightarrow B$ is decomposed as follows:
\begin{equation}
  \xymatrix{ & A \ar[dl]_{f} \ar@{-->}[r] & A' \ar[dr]^{f'} & \\
  B \ar@{-->}[rrr] & &  & B,}
\end{equation}
  where $A\dashrightarrow A'$ is one flop and
  both $f$ and $f'$ are the blow-ups 
 along sextic normal rational curves $C$ and $C'$ on $B$, respectively. 
We remark that 
this diagram already appeared 
in \cite[the proof of Lemma 4.0.4]{TZb} to show that the degree of $q_{\sS^{+}_{4}}$ is at most two.
Indeed, the rational deck transformation
$J'\colon \widetilde{\sS}^+_{4}\dashrightarrow \widetilde{\sS}^+_4$ of the map
$q_{\sS^{+}_{4}}$ is induced from the correspondence 
between $C$ and $C'$
(if the pairs $(B,C)$ and $(B,C')$ were isomorphic up to the $G$-action,
then $q_{\sS^{+}_{4}}$ would be birational, hence $J'$ would be the identity.
We show, however, this is not the case).

It is easy to see that 
the morphism 
$\Theta\colon \widetilde{U}_0\to (\mP^2)^{6}/\mathfrak{S}_6$ is 
$G$-equivariant.
Thus we can translate the study of the rational map 
$p_{\sS^{+}_{4}}\colon \sH\dashrightarrow \widetilde{S}^+_4$  
into the study of the quotient of $(\mP^2)^{6}/\mathfrak{S}_6$ by $G$.
We carefully choose a $G$-invariant open subset of 
$(\mP^2)^{6}/\mathfrak{S}_6$ such that its quotient by $G$ exists
and an involution $J$ is induced on the quotient from $J'$ through $\Theta$
(see the subsection \ref{subsection:bir} in detail).
Only in this introduction, we denote by $\sM$ this quotient.
The variety $\sM/J$ is birational to ${\sS}^{+}_{4}$.
We can study $\sM/J$ relating it with the classically studied
GIT quotient 
$\sY:=(\mP^2)^{6}/\!/\PGL_3$,
which is a compactification of 
the moduli space of ordered six distinct points on $\mP^2$.
First, $J$ has a nice interpretation.
It is classically known that $\sY$ 
has an involution called the association map.
This involution descends to an involution $j$ on $\sY/\mathfrak{S}_6$.
In Lemma \ref{nori}, we show that $J$ is nothing but a lifting of $j$.
Second, 
the $G$-action on $\sH^B_1\simeq \mP^2$ realizes $G$ as a closed subgroup
of the automorphism group $\PGL_3$ of $\mP^2$.
Indeed, $G$ is the subgroup of $\PGL_3$ consisting of elements
which preserve one fixed conic on $\mP^2$,
hence $\PGL_3/G\simeq \mP^5$ (Proposition \ref{closure}).
This implies that $\sM/J$ is birationally a $\mP^5$-bundle
over $(\sY/\mathfrak{S}_6)/j$
(Lemma \ref{serre} and the beginning of the proof of Theorem \ref{minimalaim}).

It is known that 
$(\sY/\mathfrak{S}_6)/j$ is rational.
This is a classical result due to 
A.~Coble, which easily follows from \cite[p.19 and 37]{DO}.
Therefore, to obtain the rationality of $\sM/J$,
we have only to show that $\sM/J$ is 
birationally a {\em locally trivial} $\mP^5$-bundle over $(\sY/\mathfrak{S}_6)/j$.
For this, we look for a sub $\mP^4$-bundle $\sD$ of  
$\sM\dashrightarrow \sY/\mathfrak{S}_6$ which is invariant by $J$.
Then $\sD$ descend to a sub $\mP^4$-bundle of  
$\sM/J\dashrightarrow (\sY/\mathfrak{S}_6)/j$ 
and the local triviality of $\sM/J\dashrightarrow (\sY/\mathfrak{S}_6)/j$ 
follows.
To find the sub $\mP^4$-bundle $\sD$,
we go back from $\sM$ to $\widetilde{\sS}^{+}_{4}$,
and we find the corresponding divisor on $\widetilde{\sS}^{+}_{4}$,
which is defined by
the class of sextic rational curves
such that two of their $6$ bi-secant lines intersect
(see Lemmas \ref{intersecting}--\ref{projectivebundle}). 
Now we have finished explanations of an outline of our proof of 
the rationality of $\sS^+_4$.

Finally, 
we would like to emphasize using geometries of $B$ is
natural and appropriate for the study of $\sS^+_4$. For,
the birational $\mP^5$-bundle structure on $\sS^+_4$ as above,
which is indispensable to show its rationality,
comes from the fact that the automorphism group of 
$B$ is isomorphic to $\PGL_2$.
Moreover, the Hilbert scheme $\sH$ of rational curves
of degree $6$ on $B$ ties 
$\sS^+_4$ and the moduli space of six points on $\mP^2$ modulo the $G$-action, and
the classically known association map has a good
interpretation by the birational selfmap $B\dashrightarrow B$.

\begin{ackn}
We are grateful to Professor K.~Takeuchi for
showing us his private big table of two ray games of
weak Fano $3$-folds.
Actually, he conjectured the existence of the
diagram (\ref{eq:(A)}) as in Proposition \ref{prop:tworay}.
We learned the rationality of
$(\sY/\mathfrak{S}_6)/j$ 
by private communications with Professor I.~Dolgachev,
to whom we are also grateful.

This joint work was supported
with Grant-in-Aid for Young Scientists (A).

\end{ackn}


\section{Rational curves on the quintic del Pezzo $3$-fold}
\label{section ben noto}

\subsection{Quintic del Pezzo $3$-fold $B$}
\label{quasinoto}~

Let $B\subset \mP^6$ be the smooth quintic del Pezzo $3$-fold.
$B$ is known to be unique and be realized as a linear section of
$G(2,5)$.
There are several other characterizations of $B$.
Here we give one of them, which is suitable for our purpose.

Let $\{F_2=0\}\subset \mP^2$ be
a smooth conic. Set
\[
\VSP(F_2,3)^o:=
\{(H_1, H_2, H_3)\mid H_1^2+H_2^2+H_3^2=F_2\}
\subset \Hilb^3 \check{\mP}^2,
\]
where $\check{\mP}^2$ is the dual plane to $\mP^2$,
thus linear forms $H_i$ ($i=1,2,3$) can be considered as
points in $\check{\mP}^2$. 
Mukai showed in \cite{Mu2} that
$B$ is isomorphic to 
the closed subset 
$\VSP(F_2,3):=\overline{\VSP(F_2,3)^o}\subset \Hilb^3 \check{\mP}^2$.
The variety $\VSP(F_2,3)$ has the natural action of
the subgroup $G$ of the automorphism group $\PGL_3$ of 
${\mP}^2$
consisting of elements which preserve $\{F_2=0\}$.
The group $G$ is isomorphic to $\PGL_2$.
By definition of $\VSP(F_2,3)^o$, it is easy to see that
$G$ acts on $\VSP(F_2,3)^o$ transitively.
Thus $B$ is a quasi-homogeneous $G$-variety.

\subsection{Lines on $B$}
\label{subsection:line}~

We summerize the known results about lines on $B$.

The dual plane $\check{\mP}^2$ as above
can be identified with the Hilbert scheme $\sH^B_1$
of lines on $B$, and,
for a point $b:=(H_1,H_2, H_3)\in \VSP(F_2,3)^o
\subset B$, the points
$H_i\in \check{\mP}^2$ ($i=1,2,3$) represent three lines through $b$.
By definition of $\VSP(F_2,3)^o$ and 
transitivity of the action of $G$ on $\VSP(F_2,3)^o$, 
it is easy to show the following claim:
\begin{cla}
\label{cla:doubly}
$G$ acts doubly transitively 
on the set of pairs of intersecting lines
whose intersection points are contained in 
$\VSP(F_2,3)^o$.
\end{cla}

Let $\check{F}_2$ be the dual quadratic form to $F_2$ and
$\Omega:=\{\check{F}_2=0\}$ is the associated conic in 
$\check{\mP}^2$.
Let $l$ be a line on $B$.
If $l\in \check{\mP}^2- \Omega$,
then $\sN_{l/B}=\sO_{l}\oplus\sO_{l}$.
If $l\in \Omega$, then
$\sN_{l/B}\simeq 
\sO_{\mP^1}(-1)\oplus \sO_{\mP^1}(1).$
Lines parameterized by $\Omega$ are called
{\em{special lines}}.
Let $\widetilde{\Omega}$ be the symmetric bi-linear form
associated to $\Omega$.
Then 
two lines $l$ and $m$ on $B$ intersect if and only if 
$\widetilde{\Omega}(l,m)=0$, 
where $l,m\in \sH^B_1$ are the points
corresponding to $l$ and $m$.

By the $G$-action on $B$, a line on $B$ maps to a line on $B$, hence
the $G$-action on $B$ induces
a $G$-action on $\sH^1_B$. 

\begin{prop}\label{closure}
The conic $\Omega\subset \check{\mP}^2$
is invariant under the induced action of $G$ on $\sH^B_1$.
Moreover,
this $G$ is exactly the closed subgroup of $\PGL_3$
whose elements preserve $\Omega$. 
In particular 
$\PGL_3/G \simeq \mP_*H^0(\check{\mP}^2,\sO_{\check{\mP}^2}(2))
\simeq \mP^5$.
If we take coordinates $x,y,z$ of $\sH^B_1$
such that $\Omega=\{x^2+y^2+z^2=0\}$, then
the map $\PGL_3 \to \mP^5$ is induced by
$g\in \PGL_3\mapsto {\empty^t\! g}{g}\in \mP^5$,
where we identify the vector space of symmetric matrices
with the vector space of conics on $\check{\mP}^2$.  
\end{prop}

\begin{proof}
By the $G$-action on $B$,
a special line is mapped to a special line,
thus $\Omega$ is invariant by the induced $G$-action
on $\sH_1^B$.

By \cite[p.154]{FH},
the closed subgroup of $\PGL_3$
whose elements fix $\Omega$ is isomorphic to $G$.

Now we consider the induced action of $\PGL_3$ on 
the space of conics $\mP_*H^0(\check{\mP}^2,\sO_{\check{\mP}^2}(2))$
on $\sH^B_1$.
Since $\PGL_3$ acts on $\mP_*H^0(\check{\mP}^2,\sO_{\check{\mP}^2}(2))$
transitively and the kernel of the map 
\begin{eqnarray*}
\PGL_3 &\to& \mP_*H^0(\check{\mP}^2,\sO_{\check{\mP}^2}(2))\\
g &\mapsto & g\cdot \Omega
\end{eqnarray*}
is nothing but $G$,
it holds that $\PGL_3/G \simeq \mP^5$.

It is easy to see the last assertion.
\end{proof}

Now we collect the results obtained by Furushima and Nakayama \cite{FuNa},
which is based on another characterization of $B$ by Mukai and Umemura
\cite{MU} as follows:
let $V$ be the vector space of binary sextic forms.
The group $\PGL_2$ acts on $V$ by the law
$g\cdot f_6(x,y)=f_6(ax+by,cx+dy)$,
where $f_6$ is a binary sextic form with variable $x$ and $y$, and
$g=\begin{pmatrix} a & b \\ c & d\\ \end{pmatrix}\in \PGL_2$.
Then $B$ is isomorphic to the closure of the $\PGL_2$-orbit 
$\PGL_2[xy(x^4-y^4)]$ in $\mP_*V$. 

Let 
$\pi\colon \mP\to \sH_{1}^{B}$ be the universal
family of lines on $B$. 
Denote by $\varphi\colon \mP\to B$ the natural projection. 
As we mentioned above,
$\varphi$ is a finite morphism of degree three
(see also \cite[Lemma 2.3 (1)]{FuNa}). 

\begin{nota}\label{nota}
For an irreducible curve $C$ on $B$, 
denote by $M(C)$ the locus $\subset \sH^B_1$ of lines intersecting $C$,
namely, $M(C):=\pi(\varphi^{-1}(C))$ with reduced structure. 
Since $\varphi$ is flat, $\varphi^{-1}(C)$ is purely one-dimensional.
If $\deg C\geq 2$, 
then $\varphi^{-1}(C)$ does not contain a fiber of 
$\pi$, thus $M(C)$ is a curve.
See Proposition \ref{prop:FN} for the description of $M(C)$
in case $C$ is a line.  
\end{nota}

\begin{prop} 
\label{prop:FN}
It holds$:$ 
\begin{enumerate}[$(1)$]
 \item
the union of special lines is
the branched locus $B_{\varphi}$ of 
 $\varphi\colon \mP\to B$.
$B_{\varphi}$ has the following properties:
\begin{enumerate}[$({1}\text{-}1)$]
\item $B_{\varphi}\in |-K_{B}|$, 
\item $\varphi^*B_{\varphi}=R_1+2R_2$,
where $R_1\simeq R_2\simeq \mP^1\times \mP^1$, and
$\varphi\colon R_1\to
B_{\varphi}$ and 
$\varphi\colon R_2\to
B_{\varphi}$ are injective, and
\item
the pull-back of a hyperplane section of
$B$ to $R_1$ is a divisor of type $(1,5)$,
\end{enumerate}
\item the image of $R_2$ by $\pi\colon \mP\to \sH_{1}^{B}$
is the conic $\Omega$, 
\item
if $l$ is a special line,
then $M(l)$ is the tangent line to $\Omega$ at $l$.
If $l$ is not a special line,
then $\varphi^{-1}(l)$ is the disjoint union of 
the fiber of $\pi$ corresponding to $l$,
and the smooth rational curve dominating a line on $\sH^B_1$.
In particular,
$M(l)$ is the disjoint union of a line and the point $l$.

{\em{By abuse of notation}}, we denote by $M(l)$ the one-dimensional part of
$M(l)$ for any line $l$. 
Vice-versa,
any line in $\mathcal{H}_{1}^{B}$ is of the form $M(l)$ for some line $l$,
and
\item 
the locus swept by lines intersecting $l$ is a hyperplane section $T_{l}$
  of $B$ whose singular locus is $l$. For every point $b$ of $T_l- l$,
  there exists exactly one line which belongs to $M(l)$ 
  and passes through $b$.
\end{enumerate}
\end{prop}
\begin{proof} See \cite{FuNa} and \cite[\S 1]{IlievB5}. 
\end{proof}

By the proof of \cite{FuNa} we see that $B$ is decomposed into three
$G$-orbits as follows:
\[
B=(B- B_{\varphi})\cup 
(B_{\varphi}- C_{\varphi})\cup C_{\varphi},
\] where 
$C_{\varphi}$ is a smooth rational normal sextic and 
if $b\in B- B_{\varphi}$ exactly three
distinct lines pass through it, if $b\in
(B_{\varphi}- 
C_{\varphi})$ exactly two distinct lines pass through it,
one of them is special,
and finally $C_{\varphi}$ is the unique closed $G$-orbit and
is the loci of $b\in B$ 
through which it passes only one line, which is special.
Moreover, 
$B_{\varphi}- C_{\varphi}=\PGL_2[xy^5]$ and
$C_{\varphi}=\PGL_2[y^6]$.

It also holds that $\VSP^o(F_2,3)=B- B_{\varphi}$.

\subsection{Rational curves on $B_{}$ of degree $d\leq 6$}~

In the rest of the section 2, 
We mainly review some of our results proved 
in \cite{TZ} or in \cite{TZb}. We point out the readers that 
for a general understanding of the content of this paper they only
need to remind themselves only the statements and the definitions 
contained in this section. Moreover many of these preliminary results 
should be of easy geometrical intelligibility.   

We denote by $\sH^B_d$ the union of irreducible components
of the Hilbert scheme whose general points
parameterize normal rational curves on $B$ of degree $d\leq 6$.
\begin{prop}
\label{prop:irred}
$\sH^B_d$ is an irreducible variety of dimension $2d$.
Moreover,
a general normal rational curve $C_d\in \sH^B_d$ 
is obtained as a smoothing of the union of
a general normal rational curve $C_{d-1}\in \sH^B_{d-1}$  
and a general line $l$ on $B$ intersecting $C_{d-1}$.
\end{prop}

\begin{proof}
See \cite[Proposition 2.5.2]{TZ}.
To show this fact,
we use the irreducibility of the Hilbert scheme of 
smooth rational curves on $G(2,5)$ of degree $d$ (see \cite{Per}).
\end{proof}

We investigate $\sH:=\sH^B_6$ a bit more.

\begin{prop}\label{gradoduno6} 
In the Hilbert scheme of curves of degree $6$,
the locus of $C$ satisfying one of the
    following is a divisor of $\sH:$
    \begin{enumerate}[$(1)$]
    \item
    $C$ is the union of a general quintic normal rational curve $C_5$
    and a general line $l$ intersecting $C_5$, or
    \item
    $C$ is a general sextic rational curve contained in a 
    general hyperplane section of $B$.
    \end{enumerate} 
\end{prop}

\begin{proof} 

If $C$ satisfies (1), then $C$ has one parameter smoothing to
a sextic normal rational curve by \cite[the proof of Proposition 2.3.2]{TZ}, and, conversely,
a general sextic normal rational curve is obtained in this way 
by Proposition \ref{prop:irred}.
Thus $C\in \sH$.
By Proposition \ref{prop:irred},
$\sH^B_5$ is irreducible and is of dimension $10$.
Moreover, since $M(C_5)$ is a curve for a $C_5$,
such $C$'s form a divisor of $\sH$.

Assume $C$ satisfies (2).
Let $H$ be the hyperplane section of $B$ containing $C$.
Then, by $-K_H\cdot C=6$, it holds that $(C^2)_H=4$.
By $H\cdot C=6$, it holds that $\sN_{C/B}\simeq 
\sO_{\mP^1}(4)\oplus \sO_{\mP^1}(6)$ or
$\sO_{\mP^1}(5)\oplus \sO_{\mP^1}(5)$.
In any case, the Hilbert scheme is smooth at $C$, and
is $12$-dimensional at $C$.
On the other hand, 
sextic rational curves satisfying (2) form
a $11$-dimensional family. Indeed, 
once we fix a smooth hyperplane section,
the family of smooth sextic rational curve on it is
$5$-dimensional, and hyperplane sections of $B$ move in a 
$6$-dimensional family.
Thus $C\in \sH$.
\end{proof}

Now we recall some other results on 
intersections of lines with rational curves of
degree $d\leq 6$.

\begin{prop}
\label{prop:Cd1}
A general element $C\in\sH^{B}_{d}$
satisfies the following conditions$:$
\begin{enumerate}[$(1)$]
\item
there exist no $k$-secant lines of $C$ on $B$ with $k\geq 3$,
\item there exist at most finitely many bi-secant lines of $C$ on
$B$, and any of them intersects $C$ simply, and
\item bi-secant lines of $C$ are mutually disjoint.
\end{enumerate}
\end{prop}
\begin{proof} See \cite[Proposition 2.4.1]{TZ}.
\end{proof}

We describe some more relations of $C$ with lines on $B$
which can be translated into the geometry of $\sH^B_1$.
More explicitly, we prove that $M:=M(C)$ is sufficiently general
if $C$ is general,
where the readers have to remind the notations given
in \ref{nota}.

We denote by $\beta_i$ ($1\leq i\leq s$) the bi-secant lines of 
a general $C\in \sH^B_d$.

\begin{prop}
\label{prop:Cd}
A general element $C\in\sH^{B}_{d}$
satisfies the following conditions$:$

\begin{enumerate}[$(1)$]
\item
$C$ intersects $B_{\varphi}$ simply,

\item
$M:=M(C)$ is an irreducible curve of degree $d$
with only simple nodes
$($if $d=1$, then 
we remind the readers that, in Proposition $\ref{prop:FN}$ $(3)$,
we abuse the notation by denoting the one-dimensional part of
$\pi(\varphi^{-1}(C))$ by
$M(C)$$)$,
\item
for a general line $l$ intersecting $C$,
$M\cup M(l)$ has only simple nodes as its singularities, 
\item
$M\cup M(\beta_i)$ has only simple nodes as its singularities, and
\item
for a general line $\alpha$ intersecting $\beta_i$,
$M\cup M(\alpha)$ has only simple nodes as its singularities. 
\end{enumerate}
\end{prop}
\begin{proof} See \cite[Proposition 2.4.4]{TZ}.
\end{proof}

\begin{cor}
\label{extraline}
For a general $C\in \sH:=\sH^B_6$,
there are two lines $\alpha_{i1}$ and $\alpha_{i2}$ intersecting 
both $C$ and $\beta_i$ outside $C\cap \beta_i$
$(1\leq i\leq s)$. Moreover 
if $i\not =k$, then $\alpha_{ij}$ is disjoint from $\beta_k$, where 
$i,k=1,\dots,s$ and $j=1,2$.
\end{cor}

\begin{proof}
This immediately follows from Proposition \ref{prop:Cd} (2) and (4).
\end{proof}

\subsection{Curve $\sH_1$ parameterizing lines 
on the blow-up $A$ of $B$ along a sextic rational curve}
\label{subsection:H1}~

Though the argument in this subsection works also for other degrees $d$,
we specialize to the degree $6$ case.
For readers' covenience, we repeat almost all the proofs.

We set $\sH:=\sH^B_6$ for simplicity of notation
as in the introduction.

\subsubsection{Construction of $\sH_1$}~

For a general $C\in\sH$, we set 
\[
\text{$\sH_1:=\varphi^{-1} C\subset \mP$ and $M:=M(C)$.}
\]
   
\begin{prop}\label{primaC}
$\sH_1$ is a smooth non-hyperelliptic trigonal curve of genus
$4$.
\end{prop}

\begin{proof}  
By Propositions \ref{prop:FN} (1) and \ref{prop:Cd} (1),
it holds that $\sH_1$ is smooth and
the ramification for $\sH_1\to C$ is simple.
Since $B_{\varphi}\in |-K_{B}|$,
we can compute $g(\sH_1)$ by the Hurwitz formula:
\[
  \text{$2g(\sH_1)-2=3\times (-2)+6\times 2$, equivalently,  
  $g(\sH_1)=4$}.      
\]
\end{proof}

\begin{cor}
    \label{orabisecanti}
The number of nodes of $M$ is $6$, whence
$C$ has $6$ bi-secant lines on $B$.
\end{cor} 

\begin{proof}    
Note that a bi-secant line of $C$ corresponds to a node of $M$. 
Thus, by Proposition \ref{prop:Cd1} (2),
the morphism $\pi_{|\sH_1}\colon\sH_1\rightarrow M$ is birational.
By Propositions \ref{prop:Cd} (2) and \ref{primaC},
it holds that $p_a(M)=\frac{(d-1)(d-2)}{2}=10$
and the number of nodes of $M$ is $10-g(\sH_1)=6$. 
\end{proof}




\subsubsection{Lines on the blow-up $A$}
\label{subsubsection:lineA}~

For a general $C\in \sH$,
we take the blow-up $f\colon A\to B$ along $C$.
Let $E$ be the $f$-exceptional divisor.
We need to study the families of curves on $A$ of degree one  
with respect to the anti-canonical sheaf of ${A}$ to give another
useful interpretation of the curve $\sH_{1}$. 

\begin{nota}
~
For $i=1,\ldots, s$ and $j=1,2$, we set
\begin{enumerate}[$(1)$] 
\item
 $\{p_{i1}, p_{i2}\}=C\cap \beta_{i}\subset B$,
\item
$\zeta_{ij}=f^{-1}(p_{ij})\subset E\subset A$, and
\item
$p'_{ij}=\beta'_{i}\cap \zeta_{ij}$.
\end{enumerate}
\end{nota}

\begin{defn}
    \label{linea}
    We say that a connected curve $l\subset {A}$ 
    is a {\em{line}} on ${A}$ if
$-K_{{A}} \cdot l=1$ and
${E}\cdot l=1$.
\end{defn}

We point out that since $-K_A=f^*(-K_{B})-E$ and $E\cdot l=1$ then
$f(l)$ is a line on $B$ {\it{intersecting}} $C$.
Based on this, we can classify lines on $A$ as follows:
\begin{prop}\label{lineeA}
A line $l$ on $A$ is one of the following curves on $A:$
\renewcommand{\labelenumi}{\textup{(\roman{enumi})}}
\begin{enumerate}
\item
the strict transform of a uni-secant line of $C$ on $B$,
or
\item
the union $l_{ij}=\beta'_{i}\cup \zeta_{ij}$ $(i=1,\ldots, s$,
$j=1,2)$.

\end{enumerate}
In particular $l$ is reduced and $p_a(l)=0$.
\end{prop}
\begin{proof} This follows from easy computations on the Chow ring of $A$.
\end{proof}

\begin{prop}\label{parameter} 
    The curve $\sH_1\subset\mP$ is the Hilbert scheme of
     the lines of $A$.
\end{prop}

\begin{proof} 
We only show that $\sH_1$ parameterizes lines on $A$.
See \cite[Corollary 4.1.8]{TZ} for a rigorous proof.

By definition of $\sH_1$,
we have
$\sH_1=\{(l,t)\mid l\in M, t \in C\cap l\} \subset M\times C$,
namely, $\sH_1$ parameterizes the pairs of a line $l$ and a point $t$ in 
$C\cap l$.
In \cite{TZ}, these pairs are called marked lines.
It is easy to see that there is one to one correspondence
between marked lines and lines on $A$.
Indeed, let $m$ be a line on $A$.
The line $m$ satisfies (1) or (2) of Proposition \ref{lineeA}.
If $m$ satisfies (1), then the image $f(m)$ of $m$ on $B$ is a uni-secant line,
thus a marked line $(f(m), C\cap m)$ is uniquely determined from $m$.
If $m=\beta'_{i}\cup \zeta_{ij}$,
then we assign the marked line $(\beta_i, p_{i3-j})$ to $m$.   
Therefore $\sH_1$ parameterizes lines on $A$.
\end{proof}

By the proof of Proposition \ref{parameter}, we have the following:

\begin{cor}
\label{cor:node}
$\pi_{|\sH_1}^{-1} (\beta_i)=
\{l_{i1}, l_{i2}\}$, where $\beta_i\in M$.
\end{cor}

\subsubsection{The theta characteristic on $\sH_1$}
\label{subsubsection:theta}~

Via the interpretation of 
$\sH_{1}$ recalled in subsection \ref{subsubsection:lineA}, 
we defined the following
incidence correspondence in \cite[Section 3.1]{TZb}:
\begin{equation}\label{correspondencebis}
  I:=
  \{(l_1, l_2)\in \mathcal{H}_1\times \mathcal{H}_1 \mid
\text{$l_1 \not =l_{2}$ and $l_1\cap l_{2}\neq\emptyset$ }\}.
\end{equation}

We denote by $\delta$ the $g^1_3$ on $\sH_1$ which defines 
$\varphi_{|\sH_1}\colon \sH_1\rightarrow C$.
Let $l$, $l'$ and $l''$ be three lines on $A$ such that 
$l+l'+l''\sim \delta$. Then the images of 
$l$, ${l}'$
and ${l}''$ are lines on $B$ through one point of $C$.
Set 
\begin{equation}\label{queen}
\theta:=(\pi_{|\sH_1})^*\sO_M(1)-\delta,
\end{equation}
\noindent where $\pi\colon\mP\rightarrow \sH^{B}_{1}=\mP^{2}$ is
the natural projection of the universal family and $M= \pi (\sH_{1})$. 
Note that $\deg \theta=3$.

\begin{prop}\label{prop:sopratheta} 
The class of 
    $\theta$ is an ineffective theta characteristic and
    $I=I_{\theta}$, where, by definition, $(x,y)\in I_{\theta}$ if and only if
$y$
    belongs to the support of the unique effective divisor of $\mid
    \theta+ x\mid$.
\end{prop}
\begin{proof}
See \cite[Proposition 3.1.2]{TZb}.    
\end{proof}


\section{Birational selfmap of $B$}~
We need to refine our understanding of the geometry of the blow-up $f\colon
A\rightarrow B$ along a general $C\in \sH$. 
The main result of this subsection is Proposition \ref{prop:tworay},
in which we construct a birational selfmap $B\dashrightarrow B$
from $f$ and describe it. This is a technical core of the proof of the
main theorem. We recommend the readers to understand only the statement
of this result first and go to the proof of the main theorem in the section 4.

For readers' convenience,
we give the definition and basic properties of flops (the subsection \ref{subsection:two-ray}), and
descriptions of auxiliary birational maps which originate from $B$
(the subsection \ref{subsection:aux}).


\subsection{Smooth $3$-fold flops}
\label{subsection:two-ray}~

\begin{defn}
\label{defn:flop}
Let $A$ be a smooth $3$-fold.
A projective morphism $g\colon A\to \overline{A}$ 
is called a {\em flopping contraction}
if $g$ is isomorphic outside the union $\gamma$ of
a finite number of 
curves (actually $\gamma$ is a tree of smooth rational curves)
and any irreducible component of $\gamma$ is numerically trivial
for $K_A$.
An irreducible component of $\gamma$ is called a {\em flopping curve}. 
If there exists a divisor $D$ numerically negative
for any irreducible component of $\gamma$,
then $g$ is called a {\em $D$-flopping contraction}.
It is well-known that,
for a $D$-flopping contraction $g$,
there exists a unique
projective morphism $g'\colon A'\to \overline{A}$
such that 
\begin{itemize}
\item
$g'$ is isomorphic outside the union $\gamma'$ of
a finite number of curves and 
any irreducible component of $\gamma'$ is numerically trivial
for $K_{A'}$,
\item
the map ${g'}^{-1}\circ g\colon A\dashrightarrow A'$ 
gives an isomorphism
between $A- \gamma$ and $A'- \gamma'$, and
\item
the strict transform $D'$ on $A'$ of $D$  
is numerically {\em positive} for any irreducible component of $\gamma'$
\end{itemize}
(see \cite{flop}).
The map ${g'}^{-1}\circ g\colon A\dashrightarrow A'$ is 
called the {\em $D$-flop} for $g$ and
the morphism $g'$ is called the {\em $D$-flopped contraction}.
An irreducible component of $\gamma'$ is called a {\em flopped curve}. 

In case where 
$\rho(A/\overline{A})=1$ (for example, $\gamma$ is irreducible),
then the $D$-flop is independent of $D$ and we say simply
$A\dashrightarrow A'$ is the flop, $g'$ is the flopped contraction, etc.
\end{defn}

In Proposition \ref{prop:symm},
we summerize basic properties of flops, 
for which it is easy to find references in the literatures:

\begin{prop}
\label{prop:symm}
Let $A$ be a smooth $3$-fold and $D$ a divisor on $A$.
Let $g\colon A\to \overline{A}$ be a $D$-flopping contraction
and $\gamma$ the union of all the flopping curves.
Let $A\dashrightarrow A'$ the $D$-flop and
$g\colon A'\to \overline{A}$ the $D$-flopped contraction.
Denote by $\gamma'$ the union of all the $D$-flopped curves.
Then 
\begin{enumerate}[$(1)$]
\item
$A'$ is smooth, 
\item
$g$ and $g'$ is isomorphic analytically near $\gamma$ and $\gamma'$.
In particular, the numbers of irreducible components of $\gamma$ and $\gamma'$
are equal, and
\item
if $\rho(A/\overline{A})=1$, then
$G\cdot \gamma=-G'\cdot \gamma'$,
where $G$ is a divisor on $A$ and $G'$ is the strict transform on $A'$ of $G$.
\end{enumerate}  
\end{prop}

\begin{proof}
See \cite{flop}.
\end{proof}

\begin{expl}[Atiyah's flop]
\label{expl:Atiyah}
Here we describe the simplest flopping contraction.
Actually, in the sequel, we mainly need only (composites of) 
flopping contractions of this type.

Let $g\colon A\to \overline{A}$ be
a projective morphism whose
exceptional curve $\gamma$ is a smooth irreducible rational
curve with $\sN_{\gamma/A}\simeq \sO_{\mP^1}(-1)\oplus
\sO_{\mP^1}(-1)$. 
It is easy to check that $g$ is a flopping contraction.
We can construct the flop $A\dashrightarrow A'$
as follows: 
let $p\colon \widehat{A}\to A$ be the blow-up of $A$ along $\gamma$
and $E$ the $p$-exceptional divisor.
Since $\sN_{\gamma/A}\simeq \sO_{\mP^1}(-1)\oplus
\sO_{\mP^1}(-1)$, it holds that $E\simeq \mP^1\times \mP^1$.
There exists a morphism $q\colon \widehat{A}\to A'$
which is isomorphic outside $E$ and
$q_{|E}$ is the natural projection 
$E\simeq \mP^1\times \mP^1\to \mP^1$ different from
$E\to \gamma$.
It is easy to check that
there exists a projective morphism $g'\colon A'\to \overline{A}$
which is isomorphic outside $\gamma':=q(E)$ and
$q\circ p^{-1}\colon A\dashrightarrow A'$
is the flop.
The flop $A\dashrightarrow A'$ is called {\em Atiyah's flop}.
\end{expl}

The following two results, Propositions \ref{flop} and \ref{prop:Atiyah},
describe changes of intersection numbers by a flop.
They are well-known for the experts but are not explicitly stated 
in the literatures. Therefore we decided to write their proofs in full details. 

\begin{prop}
\label{flop}
Let $A$ be a smooth $3$-fold
and $g\colon A\to \overline{A}$ be
a flopping contraction with $\rho(A/\overline{A})=1$.
Denote by $\gamma$ the union of all the $g$-exceptional curves.
On $A$, take a divisor $N$ and
an irreducible projective curve $\delta\not \subset \gamma$.
Let $A\dashrightarrow A'$ be the flop, and 
$N'$ and $\delta'$ the strict transforms on $A'$ of $N$ and
$\delta$ respectively. 
It holds:
\begin{enumerate}[$(1)$]
\item
If $N\cdot \gamma =0$, then 
$N^3={N'}^3$ and 
$N\cdot \delta = N'\cdot \delta'$.
\item
If $N\cdot \gamma >0$,  
then $N^3>{N'}^3$
and $N\cdot \delta \leq N'\cdot \delta'$.
\item
If $N\cdot \gamma <0$,  
then $N^3<{N'}^3$
and $N\cdot \delta \geq N'\cdot \delta'$.
\end{enumerate}
\end{prop}

\begin{proof}
By Proposition \ref{prop:symm}, 
the inverse $A'\dashrightarrow A$ of $A\dashrightarrow A'$
is also the flop
for the flopping contraction $g'\colon A'\to \overline{A}$,
thus
we may assume that $N\cdot \gamma\geq 0$
by interchanging the roles of $A$ and $A'$.

First we verify the inequality between $N^3$ and ${N'}^3$.
We learned the proof by \cite[Corollary 9.3]{SShB2},
which originated from Mori.
We write the proof for readers' convenience.

Since we assume $N$ is $g$-gef, 
$\Bs |m(N+H)|=\emptyset$
by Kawamata-Shokurov's base point free theorem (\cite[Theorem 3-1-1]{KMM}),
where $m\gg 0$ and
$H$ is the pull-back of a sufficiently ample divisor on $\overline{A}$. 
Take $H_1,H_2,H_3\in |H|$ and
$N_1,N_2,N_3 \in |m(N+H)|$ such that
$H_i$ are disjoint from $\gamma$ and
$N_i$ do not intersect each other on $\gamma$.
For any divisor $L$ on $A$,
we denote by $L'$ its strict transform on $A'$.
It holds that $L_1\cdot L_2\cdot H_i=L'_1\cdot L'_2\cdot H'_i$
for any divisors $L_1$ and $L_2$ on $A$
since $H_i\cap \gamma =\emptyset$ and $A\dashrightarrow A'$
is isomorphic outside $\gamma$.
Then we have
\begin{eqnarray*}
m^3(N^3-{N'}^3)= m^3\{(N+H)^3-(N'+H')^3\}
=N_1\cdot N_2\cdot N_3-N'_1\cdot N'_2\cdot N'_3=\\
(N_1\cdot N_2\cdot N_3)_{\gamma}-
(N'_1\cdot N'_2\cdot N'_3)_{\gamma'}=
-(N'_1\cdot N'_2\cdot N'_3)_{\gamma'}
\end{eqnarray*}
If $N\cdot \gamma=0$, then we may assume that 
$N_i\cap \gamma=\emptyset$, hence $N'_i\cap \gamma'=\emptyset$.
This implies that
$m^3(N^3-{N'}^3)=-(N'_1\cdot N'_2\cdot N'_3)_{\gamma'}=0$.
If $N\cdot \gamma>0$, then
it holds that $N'\cdot \gamma'<0$ by Proposition \ref{prop:symm} (3).
Thus $\Supp (N'_2\cap N'_3)=\gamma'$ and
$m^3(N^3-{N'}^3)=-(N'_1\cdot N'_2\cdot N'_3)_{\gamma'}=
-N'_1\cdot (N'_2 \cdot N'_3)_{\gamma'}>0$.

Second we verify the inequality between 
$N\cdot \delta$ and ${N'}\cdot \delta'$.
Take the following diagram:
\begin{equation}
\label{eq:line}
\xymatrix{ & \widehat{A} \ar[dl]_{p} \ar[dr]^{q}\\
 A &  & A',}
\end{equation}
where $p$ and $q$ are resolutions of $A$ and $A'$
respectively.
By definition of the flop,
$A\dashrightarrow A'$ is isomorphic outside
$\gamma$ and $\gamma'$.
Therefore
we may assume that $p$ (resp.~$q$) is
isomorphic outside $\gamma$ (resp.~$\gamma'$).
We can write
$q^*N'=p^*N+R$,
where $R$ is a $p$-exceptional, hence is also a $q$-exceptional divisor.
The divisor $p^*N$ is $q$-nef since we assume that $N$ is $g$-nef.
By the negativity lemma (cf. \cite[Lemma 2.19]{FA}), 
it holds $R\geq 0$.
The inequality $N\cdot \delta \leq N'\cdot \delta'$
follows from this fact.

Assume that $N\cdot \gamma=0$. Then, by Proposition \ref{prop:symm},
it holds that $N'\cdot \gamma'=0$.
Therefore we can interchange the role of $A$ and $A'$ and
we have also $R\leq 0$ by
applying the negativity lemma to
$p^*N=q^*N'-R$.
Consequently, we have that $p^*N=q^*N'$ and
then the equality $N\cdot \delta = N'\cdot \delta'$.

\end{proof}

Now we specialize to Atiyah's flop and refine Proposition \ref{flop}.

\begin{prop}
\label{prop:Atiyah}
Let $A$ be a smooth $3$-fold
and $g\colon A\to \overline{A}$ be
a flopping contraction 
whose exceptional curve $\gamma$ is irreducible.
Assume that $\sN_{\gamma/A}\simeq \sO_{\mP^1}(-1)\oplus \sO_{\mP^1}(-1)$.
Let $N$ be a divisor on $A$ and set $d:=N\cdot \gamma$.
Let $\delta$ be an smooth irreducible projective curve different
from $\gamma$ and set $e$ to be the set-theoretic
intersection number of $\delta$ and 
$\gamma$. 
Let $A\dashrightarrow A'$ be the flop, and 
$N'$ and $\delta'$ the strict transforms on $A'$ of $N$ and
$\delta$ respectively. 

It holds that
$N^3=(N')^3+d^3$ and
$N'\cdot \delta' \geq N\cdot \delta+de$.
Moreover, if $\gamma$ and $\delta$ intersect transversely at $e$ points,
then $N'\cdot \delta'=N\cdot \delta+de$. 
\end{prop}

\begin{proof}
Take the following diagram as in Example \ref{expl:Atiyah}:
\begin{equation}
\label{eq:line}
\xymatrix{ & \widehat{A} \ar[dl]_{p} \ar[dr]^{q}\\
 A &  & A',}
\end{equation}
where $p$ is the blow-up along $\gamma$
and $q$ is the blow-down of $p$-exceptional divisor $E\simeq \mP^1\times
\mP^1$ in the other direction.
We can write
$q^*N'=p^*N+aE$ with some $a\in \mZ$.
We show that $a=d$. Indeed,
for a fiber $\widehat{\gamma}$ of $E\to \gamma'$,
which is mapped to $\gamma$ by $p$,
it holds
\[
\text{$q^*N'\cdot \widehat{\gamma}=0$,\,
$p^*N\cdot \widehat{\gamma}=N\cdot \gamma=d$,\, and
$E\cdot \widehat{\gamma}=-1$.}
\]
Therefore we have $a=d$.

Now we prove the inequality 
$N'\cdot \delta' \geq N\cdot \delta+de$.
Let $\widehat{\delta}$ be the strict transform on $\widehat{A}$
of $\delta$.
By definition of $e$,
it holds that $E\cdot\widehat{\delta}\geq e$.
By $q^*N'=p^*N+dE$, we have
\[
N'\cdot \delta'=q^*N'\cdot \widehat{\delta}=
(p^*N+dE)\cdot \widehat{\delta}\geq 
p^*N\cdot \widehat{\delta}+de=N\cdot \delta+de.
\]
Moreover, if $\gamma$ and $\delta$ intersect transversely at $e$ points,
then it holds that $E\cdot\widehat{\delta}=e$.
Thus we have $N'\cdot \delta'=N\cdot \delta+de$.

To prove the equality $N^3=(N')^3+d^3$,
we compute $p^* N^2 q^*N'$ in two ways.
First, by applying the projection formula to $p$,
 we have $p^* N^2 q^*N'=N^3$.
Second, by the equality $p^*N=q^*N'-dE$,
we have 
\[
p^* N^2 q^*N'=(q^*N'-dE)^2q^*N'
=(q^*N')^3+d^2 E^2 q^*N'=(N')^3+d^2 N'\cdot q_*(E^2),
\]   
where it holds that $(q^*N')^2 E=(q^*N'_{|E})^2=0$
since $E$ is a $\mP^1$-bundle over a curve and
$q^*N'_{|E}$ is numerically a sum of its fibers.
Thus we have
$N^3=(N')^3+d^2N'\cdot q_*(E^2)$.   
It is easy to see that
$-q_*(E^2)=\gamma'$ as a $1$-cycle.
Therefore $N'q_*(E^2)=-N'\cdot \gamma'=N\cdot \gamma=d$ 
by Proposition \ref{prop:symm} (3).
Consequently, we have the equality $N^3=(N')^3+d^3$.   

\end{proof}

\subsection{Auxiliary birational maps originating from $B$}
\label{subsection:aux}

\begin{prop}
\label{projline}
Let $l$ be a line on $B$. Then the 
projection of $B$ from $l$ is decomposed as follows$:$
\begin{equation}
\label{eq:line}
\xymatrix{ & B_l \ar[dl]_{\pi_{1 l}} \ar[dr]^{\pi_{2 l}}\\
 B &  & Q, }
\end{equation}
where $\pi_{1l}$ is the blow-up along $l$ and 
$B \dashrightarrow Q$ is the projection from $l$ and 
$\pi_{2l}$ contracts onto a twisted cubic curve the
strict transform of the locus $T_l$ swept by the lines of $B$
touching $l$. Moreover 
\begin{equation}
\label{eq:line2}
-K_{B_l}=H_l+L_l,
\end{equation}
where $H_l$ and $L_l$ are the pull backs of general hyperplane sections of 
$B$ and $Q$ respectively.
We denote by $E_l$ the $\pi_{1l}$-exceptional divisor.
\end{prop}
\begin{proof} 
This is well-known and explicitly stated in 
\cite{Fu2} and \cite{MM1}.
See also \cite[Proposition 3.1.1]{TZ}.
\end{proof}

As an application, we show the following,
which we need in the section 4:
\begin{cor}
\label{cor:gen}
For a general $C\in \sH$, the six points 
$\beta_{1},\ldots,\beta_{6}$ on $\mP^2$ are in a general position.
\end{cor}

\begin{proof}
Assume by contradiction that
there exists a line $L$ through a set of $3$ points $\beta_{i_j}\in \mP^2$
$(1\leq j\leq 3)$.
By Proposition \ref{prop:FN} (3),
there exists a line ${l}$ on $B$ such that $M({l})=L$.
The above condition means that
$3$ bi-secant lines $\beta_{i_j}$ $(1\leq j\leq 3)$ intersect ${l}$.
Consider the successive linear projections 
$B\dashrightarrow Q\dashrightarrow \mP^2$
first from $\beta_{i_1}$
and then from the strict transform on $Q$ of $\beta_{i_2}$.
The image $\widetilde{C}$ of $C$ on $\mP^2$ is a line or a conic.
If $\widetilde{C}$ is a line, then $C$ is contained in a hyperplane section,
a contradiction.
Thus $\widetilde{C}$ is a conic and $C\dashrightarrow \widetilde{C}$ 
is an isomorphism.
However,
the images of $\beta_{i_2}$ and $\beta_{i_3}$ on $Q$ mutually intersect
since $l$ is contracted by $B\dashrightarrow Q$.
Thus the image on $Q$ of $\beta_{i_3}$ is contracted by the projection
$Q\dashrightarrow \mP^2$.
Moreover, the images of $\beta_{i_2}$ and $\beta_{i_3}$ on $Q$ 
are bi-secant lines of the image of $C$,
hence $\widetilde{C}$ must be singular at the image of $\beta_{i_3}$,
a contradiction.

In the proof of \cite[Lemma 3.1.1]{TZb},
we have shown that there are no conic through
the six points 
$\beta_{1},\ldots,\beta_{6}$
using the inductive construction of $C$. 
\end{proof}

\begin{defn}
Let $b$ be a point of $B$.
We call the rational map from $B$
defined by the linear system of hyperplane sections
singular at $b$ {\em{the double projection from $b$}}.
\end{defn}

\begin{prop}
\label{prop:proj2}
For a point $b\in B- B_{\varphi}$,
the double projection from $b$ is described as follows:
\begin{enumerate}[$(1)$]
\item
 
the target of the double projection is $\mP^2$,
and the double projection from $b$ and 
the projection $B\dashrightarrow \overline{B}_b$ from $b$ 
fit into the following diagram:
\begin{equation}
\label{eq:pt1}
\xymatrix{ & B_b \ar[dl]_{\pi_{1 b}} \ar[dr]\ar@{-->}[rr] & & 
 B'_b \ar[dl] \ar[dr]^{\pi_{2 b}} & \\
 B &  & \overline{B}_b & & \mP^2, }
\end{equation}
where $\pi_{1 b}$ is the blow-up of $B$ at $b$,
$B_b\dashrightarrow B'_b$
is the flop of the strict transforms of three lines through $b$, 
and $\pi_{2 b}\colon B'_b\rightarrow \mP^{2}$ 
is a $(\text{unique})$ $\mP^1$-bundle structure.

We denote by
$E_b$ the $\pi_{1b}$-exceptional divisor.
For simplicity of notation,
we denote the strict transforms on $B'_b$ of divisors on $B_b$
by the same notation.
\item 
\label{eq:pt2}
\[
L_b=H_b-2E_b\ \text{and}\ -K_{B'_b}=H_b+L_b,
\]
where 
$H_b$ is the strict transform 
of a general hyperplane section of $B$,
and $L_b$ is the pull back of a line on $\mP^2$,
\item
the strict transforms $l'_i$ of three lines $l_i$ through $b$ on $B_b$
have the normal bundle $\sO_{\mP^1}(-1)\oplus \sO_{\mP^1}(-1)$.
The flop $B_b\dashrightarrow B'_b$ is Atiyah's flop.
\item
a fiber of $\pi_{2 b}$ not contained in $E'_b$
is the strict transform of a conic through $b$, 
or
the strict transform of a line $\not \ni b$ intersecting a line through $b$.
\end{enumerate}
\end{prop}
\begin{proof} 
This is well-known for the experts but is not explicitly
stated in the literatures. See 
\cite[the proof of Proposition 3.2.2]{TZ}
for a sketch of its proof.
\end{proof}

\subsection{The birational selfmap of $B$}~

The following is one of consequences of generality of 
a sextic normal rational curve.
\begin{prop}
\label{-1-1}
Let $C$ be a general sextic normal rational curve on $B$
and $\beta_i$ $(1\leq i\leq 6)$ its six bi-secant lines.
Let $f\colon A\to B$ be the blow-up of $B$ along $C$
and $\beta'_i$ the strict transforms on $A$ of $\beta_i$.
Then $\sN_{\beta'_i/A}\simeq \sO_{\mP^1}(-1)\oplus \sO_{\mP^1}(-1)$.
\end{prop}

\begin{proof}
See \cite[Lemma 3.1.4]{TZ}.
\end{proof}

Now we reach the main result of the section 3.
The method of its proof we take was developed more or less 
by Takeuchi in the paper \cite{Take89}.
We write a proof in full details hoping 
that it becomes a good
introduction to the readers of a method of the explicit $3$-fold Mori theory.

\begin{prop}
\label{prop:tworay}
Let $C$ be a sextic normal rational curve on $B$
and $f\colon A\to B$ the blow-up along $C$.
There exists possibly a $5$-dimensional locus $\sS$ in $\sH$
such that
if $C\not \in \sS$, then we have the following diagram
\begin{equation}
\label{eq:(A)}
  \xymatrix{ & A \ar[dl]_{f} \ar@{-->}[r] & A' \ar[dr]^{f'} & \\
  B  &  &  & B,}
\end{equation}
  where $A\dashrightarrow A'$ is one flop and
  $f':A'\to B$ is also the blow-up 
 along a sextic normal rational curve $C'$. 
Denote by $E'$ the $f'$-exceptional divisor.
For simplicity of notation,
we denote the strict transforms on $A'$ of curves and divisors on $A$
by the same notation.
It holds
\begin{equation}
 \label{eq:(B)}
\text{$L=3H-2E$, $-2K_{A}=H+L$ and $E'=4H-3E$},  
\end{equation}   
  where $H$ $($resp. $L)$ is 
  the strict transform of a general hyperplane section of the
  $B$ on the left $($resp. right$)$ hand side.

Moreover, if $C$ is general, then
it holds:
\begin{enumerate}[$(1)$] 
\item
all the flopping curves of $A\dashrightarrow A'$ 
are the six strict transforms 
$\beta'_{1},\ldots ,\beta'_{6}$ of 
six bi-secant lines $\beta_{1},\ldots ,\beta_{6}$ of $C$, and
\item a non-trivial fiber
of $f'$ is the strict transform
of an irreducible tri-secant conic of $C$, or
a line intersecting both $C$ and 
a bi-secant line $\beta$ of $C$ outside $C\cap \beta$.
\end{enumerate}
\end{prop}

\begin{proof}
We divide the proof into several steps.
For applying the Mori theory,
the first step is
to check $A$ is a weak Fano $3$-fold, namely,
$-K_A$ is nef and big.
Then we can carry on the so-called two-ray game (a special case of
the minimal model program).
In the present case, we have more;
$|-K_A|=|2H-E|$ is base point free since
$C$ is the intersection of quadrics.
The bigness of $-K_A$ follows by the calculation:
\[
(-K_A)^3=(f^*(-K_B)-E)^3=(2H-E)^3=8H^3+6HE^2-E^3=14>0,
\]
where we use basic numerical equalities:
\begin{equation}
\label{eq:num}
\text{$H^3=5$, $H^2E=0$, $HE^2=-6$ and $E^3=-10$}.
\end{equation}

Let $g\colon A\to \overline{A}$ be the Stein factorization
of the morphism defined by $|-K_A|$.
We need to make a case division.\\
{\bf Case 1.} $g$ contracts a divisor $F$.

We show that such $C$'s satisfying
this condition form at most a $5$-dimensional family in $\sH$.

We can write $F\sim aH-bE$, where $a,b\in \mZ$.
It holds that $(-K_A)^2 F=0$.
By $-K_A=2H-E$ and (\ref{eq:num}), we have
$(-K_A)^2 F=14(a-b)=0$. Thus $F=a(H-E)$.
The image $g(F)$ of $F$ is not a point since $-K_A F^2=-4a^2\not =0$. 
For a fiber $l$ of $F\to g(F)$,
it holds that
$F\cdot l=-1$ or $-2$.
If $F\cdot l=-1$, then $a=1$ and $F\sim H-E$.
This is impossible; $|H-E|$ is empty
since $C$ is not contained in a hyperplane section,
Thus $F\cdot l=-2$ and $F\sim 2(H-E)$.
Together with the equality $-K_A\cdot l=(2H-E)\cdot l=0$,
it holds that $H\cdot l=1$ and $E\cdot l=2$, namely,
$l$ is irreducible and is the strict transform
of a bi-secant line of $C$.
Now we consider the set-up as in the subsection \ref{subsection:line}.
Let $\Lambda$ be the curve in $\sH^B_1$ parametrizing lines
which are the images on $B$ of fibers of $F\to g(F)$.
It holds that $\Lambda$ is an irreducible conic since
$f(F)\cdot m=2$ for a general line $m$ on $B$ and
there exists one image of a fiber of $F\to g(F)$ through
one point of $f(F)$.
We show that $C$ is determined from $\Lambda$.
Then we are done since 
conics in $\sH^1_B$ form a $5$-dimensional family.
Let $F'\subset \mP$ be the pull-back of $\Lambda$ 
by $\pi\colon \mP\to \sH^1_B$. 
Then $\varphi(F')=f(F)$.
Moreover, $\varphi_{|F'}\colon F'\to f(F)$ is birational
since there exists one image of a fiber of $F\to g(F)$ through
one point of $f(F)$.
The natural morphism $F\to f(F)$ is an isomorphism outside
$F\cap E$. For a fiber $\gamma$ of $E\to C$,
it holds that $F\cdot \gamma=2(H-E)\cdot \gamma=2$,
thus $f(F)$ is singular along $C$.
Therefore $C$ is determined from $\Lambda$ as the singular locus of
the image of $F'$ by $\varphi$.

From now on we assume that $C$ does not belong to
such a $5$-dimensional family. Thus we fall into the following case:\\
{\bf Case. 2.} $g$ contracts only finite number of curves.

Then $g$ is a flopping contraction. 
Moreover, it holds that $\rho(A/\overline{A})=1$ since $\rho(A)=2$.
Let $A\dashrightarrow A'$ be the flop.
Since $A'$ is rational,
$K_{A'}$ is not nef.
Therefore there exists an extremal contraction
$f'\colon A'\to B'$.
The morphism $f'$ is unique since $\rho(A')=2$.
For simplicity of notation,
we denote the strict transforms on $A'$ of curves and divisors on $A$
by the same notation.
We would like to determine the type of $f'$ as in the statement of 
this theorem.\\
{\bf Step. 1.}
Let $L:=3H-2E$. We show that $L$ is nef on $A'$
and $f'$ is the Stein factorization
of the morphism defined by some multiple of $L$.

We see that there exists no effective divisor $D\sim aH-bE$ on $A$
such that $a>0$ and $b\geq a$.
Indeed, if such a $D$ exists,
then $(-K_A)^2 D\leq 0$ by (\ref{eq:num}), hence $(-K_A)^2 D=0$ and
$D$ is the $g$-exceptional divisor since $-K_A$ is nef.
This contradicts the assumption that $g$ is a flopping contraction.
Thus any nonzero effective divisor $D\sim aH-bE$ satisfies
that $a=0$ and $b<0$, or $a>0$ and $b<a$.

We show that $|L|$ has no fixed component.
Assume by contradiction that $|L|$ has a fixed component.
If $E$ is a fix component, then $L-E\sim 3H-3E$ is effective,
a contradiction to the above consideration.
If there exists 
a fixed component $D\sim aH-bH$ with $a>0$ and $b<a$,
then 
$L-(aH-bE)=(3-a)H-(2-b)E$ is effective,
thus $3-a>0$ and $2-b<3-a$.
The inequality $b<a$ and $2-b<3-a$ has no solution, a contradiction.
Therefore $|L|$ has no fixed component.

We prove that $h^0(\sO_A(L))\geq 7$.
Consider the exact sequence
\begin{equation}
\label{eq:exact1}
0\to \sO_A(L)\to \sO_A(3H-E)\to \sO_E(3H-E)\to 0.
\end{equation}
$3H-E$ is nef since $3H-E=2H-E+H=-K_A+H$, and
$-K_A$ and $H$ are nef.
Thus, by the Kawamata-Viewheg vanishing theorem,
$h^0(\sO_A(3H-E))=\chi(\sO_A(3H-E))
=\frac{1}{12}(120H^3+49HE^2-6E^3)+\frac{1}{12} H\cdot c_2(A)+3$.
Let $H_0\in |H|$ be a general member.
By the exact sequence
\[
0\to T_{H_0}\to T_{A|H_0}\to \sO_{H_0}(H)\to 0,
\]
we can calculate
$c_2(A)\cdot H=18$.
Thus, by (\ref{eq:num}), we have 
$h^0(\sO_A(3H-E))=35$.
Now we compute 
$h^0(\sO_E(3H-E))$.
Note that $E$ is a $\mP^1$-bundle over $C\simeq \mP^1$.
Let $l$ be a fiber of $E\to C$.
Then $(3H-E)\cdot l=1$.
Thus $f_{|E*} \sO_E(3H-E)=\sO_{\mP^1}(a)\oplus \sO_{\mP^1}(b)$,
where $a+b=(3H-E)^2 E=26$ and $a,b\geq 0$ since $3H-E$ is nef.
Thus $h^0(\sO_E(3H-E))=28$.
Finally we have  
$h^0(\sO_A(L))\geq 7$ from (\ref{eq:exact1}).

Now we prove that $L$ is nef on $A'$.
Since $\rho(A')=2$,
it suffices to check that
$L$ is non-negative both for a flopped curve and
a general curve in a general fiber of $f'$. 
First we check that $L$ is positive for a flopped curve on $A'$.
Indeed, for a flopping curve $\gamma$,
it holds that $H\cdot \gamma>0$ and 
$(2H-E)\cdot \gamma=-K_A\cdot \gamma=0$.
Thus $L\cdot\gamma=(3H-2E)\cdot\gamma<0$.
Then, by Proposition 3.2,  
$L$ is positive for a flopped curve on $A'$.
Second we check $L$ is non-negative 
for a general curve in a general fiber of $f'$.
If $f'$ is of fiber type, then
curves in fibers cover $A'$ whence
$L$ is non-negative for a general curve in a general fiber of $f'$
since $|L|\not =\emptyset$.
If $f'$ is birational,
then, again, 
$L$ is non-negative for a general curve in a general fiber of $f'$
since the $f'$-exceptional divisor is not
a fixed component of $|L|$ on $A'$. 

Finally we show that $f'$ is defined by some multiple of $L$.
For this we prove the existence of 
an irreducible $k$-secant conic of $C$
with $k\geq 3$ by 
the double projection from a general point $b$ of $C$.
We may assume that $C$ is not contained in 
$B_{\varphi}$. 
Indeed,  
if $C$ is contained in $B_{\varphi}$,
then the pull-back of $C$ on $R_1$ is a divisor of type $(1,1)$
by Proposition \ref{prop:FN} (1-3) and $\deg C=6$.
Thus such $C$'s form $3$-dimensional family
(we do not prove the existence of such $C$'s).
We may assume that $C$ does not belong to this $3$-dimensional family.
Thus we may assume that $b\not \in B_{\varphi}$
and then there are three lines $l_1$, $l_2$ and $l_3$
through $b$.
We consider the double projection from $b$ and we use the notation
of Proposition \ref{prop:proj2}.
Since $C$ has only finitely many bi-secant lines,
we may assume that $l_i$ are not bi-secant lines
by generality of $b$.
Thus the strict transforms $C'$ and $l'_i$ of 
$C$ and $l_i$ are disjoint on $B_b$.
By $-K_{B_b}=\pi_{1b}^*(-K_B)-2E_b$, it holds that
$-K_{B_b}\cdot C'=10$.
Thus
it holds that $H_b\cdot C'=6$ on $B'_b$ and $-K_{B'_b}\cdot C'=10$,
where we denote by $C'$ the strict transform on $B'_b$ of $C'$
abusing the notation.
Hence $L_b\cdot C'=4$ by Proposition \ref{prop:proj2} (2)
and then the image of $C'$ on $\mP^2$ is
a line, a conic or a quartic.
This implies that $\pi_b$ has a multi-secant fiber of $C'$.
If it is the strict transform of a smooth conic $q$ through $b$,
then $q$ is a $k$-secant conic of $C$
with $k\geq 3$.
Otherwise, the fiber is the strict transform of
a bi-secant line of $C$ intersecting one of $l_i$.
We show this does not occur if $b$ is general.
If this occurs for general $b$'s,
then $C$ is contained in the locus
of lines $T_{\beta}$
intersecting one fixed bi-secant line $\beta$
since there are a finite number of bi-secant lines
of $C$.
This is a contradiction since $C$ is not contained in a hyperplane section.
Therefore there exists an irreducible $k$-secant conic of $C$
with $k\geq 3$.
 
Let $q$ be a general irreducible $k$-secant conic of $C$
with $k\geq 3$.
Then $L\cdot q=6-2k$ on $A$.
Since a flopping curve of $A\dashrightarrow A'$ intersects $L$ negatively,
we have $L\cdot q\leq 6-2k$ on $A'$
by Proposition \ref{flop} (3).
Since $L$ is nef on $A'$, 
we have $k=3$ and $L\cdot q=0$ on $A'$.
Thus $L$ is not ample.
By Kawamata-Shokurov's base point free theorem (\cite[Theorem 3-1-1]{KMM}),
some multiple of $L$ defines a morphism,
which is non-trivial since $L$ is not ample. 
The extremal contraction $f'$
is nothing but the Stein factorization of
the morphism defined by some multiple of $L$.

To determine $f'$,
we make a case division using the classification of extremal contractions from
smooth $3$-folds \cite{ray}.
Note that $L$ is the pull-back of a generator
of $\Pic B'$ since $L$ is primitive.\\
{\bf Step 2.}
We exclude the case where $f'$ is of fiber type.

Then $B'\simeq \mP^1$ or $\mP^2$.
We can derive this fact as follows:
it is well-known that $B'$ is smooth if $f'$ is of fiber type \cite{ray}.
Since $A$ is rational, $B'$ is covered by rational curves, thus is rational
since $\dim B'\leq 2$. If $\dim B'=1$, then $B'\simeq \mP^1$.
If $\dim B'=2$, then $B'\simeq \mP^2$ 
since the Picard number of $B'$ is one.
Thus $L$ is the pull-back of a point or a line respectively.
This is a contradiction since $h^0(L)\geq 7$
as in Step 1.\\
{\bf Step 3.}
Assuming that $f'$ contracts a divisor $E'$ to a curve $C'$,
we show that $f'$ is described as in the statement of the theorem.

$B'$ is a smooth Fano $3$-fold with $\rho(B')=1$.
By the classification of smooth Fano $3$-folds
with Picard number one,
we may write $-K_{B'}=a\widehat{L}$ with $a=1,2,3,4$,
where $\widehat{L}$ is the image of $L$ by $f'$.
Equivalently, it holds that ${f'}^*(-K_{B'})=aL$. 
If $a=3$, then $B'$ is the quadric $3$-fold
and if $a=4$, then $B'\simeq \mP^3$.
These cases contradict $h^0(L)\geq 7$.
If $a=1$, then $-K_{A'}={f'}^*(-K_{B'})-E'=2H-E$,
thus $E'=H-E$, a contradiction since $h^0(H-E)=0$.
If $a=2$, then, by the inequality $h^0(L)\geq 7$
and the classification of del Pezzo $3$-folds (see \cite{Fu}),
we have $h^0(\widehat{L})=7$ and ${\widehat{L}}^3=5$ on $B'$.
Thus $B'$ is also the quintic del Pezzo $3$-fold.
We can easily show that $C'$ is a sextic normal rational curve.
We check the equalities (\ref{eq:(B)}).
By definition of $L$, we have the former two equalities.  
By $-K_{A'}=2L-E'$, $-K_{A'}=2H-E$ and $L=3H-2E$,
we have the latter equality.
Assuming $C$ is general, we check the assertions (1) and (2).
Actually, it suffices to assume that
$C$ has six mutually disjoint bi-secant lines $\beta_i$ ($1\leq i\leq 6$)
(Proposition \ref{prop:Cd1} (3) and Corollary \ref{orabisecanti}),
and
$\sN_{\beta'_i/A}\simeq \sO_{\mP^1}(-1)\oplus \sO_{\mP^1}(-1)$
for the strict transform $\beta'_i$ of $\beta_i$ 
(Proposition \ref{-1-1}).
Any $\beta'_i$ is a $g$-exceptional curve.
We show that $\beta'_i$ ($1\leq i\leq 6$) are
the only $g$-exceptional curves.
Passing to the analytic category and taking the algebraization, 
we can decompose the flop $A\dashrightarrow A'$
into a sequence of flops
$A:=A_1\dashrightarrow A_2\dashrightarrow \cdots \dashrightarrow A_n=:A'$
for some $n\in \mN$,
where $A_j\dashrightarrow A_{j+1}$ is the flop of 
the strict transform of $\beta'_j$ if $1\leq j\leq 6$, or
the flop of the strict transform of an irreducible
$g$-exceptional curve different from $\beta'_i$ $(1\leq i \leq 6)$
if $6<j\leq n-1$. 
For simplicity of notation,
we denote by the same notation 
the strict transforms of $g$-exceptional curves, $L$ and $H$
on each $A_j$. 
Noting $L=3H-2E$, we can easily compute that
$L^3=-1$ on $A$.
Since $L$ on $A'$ is the pull-back of $\widehat{L}$,
we have $L^3=5$ on $A'$.
Note that the flop $A_j\dashrightarrow A_{j+1}$ $(1\leq j\leq 6)$
is Atiyah's flop.
Thus by the equality $L\cdot \beta'_i=-1$ ($1\leq i\leq 6$) on $A$,
we see that $L^3=-1+6=5$ on $A_7$
by Proposition \ref{prop:Atiyah}.
Assume by contradiction that there exists at least one $g$-exceptional curve
different from $\beta'_i$'s, namely, $n>7$.
Since  
the strict transforms of 
all the other $g$-exceptional curves are still numerically negative
for $L$ on $A_j$ ($j\geq 7$) by Proposition \ref{flop} (3),
it holds that $L^3>5$ on $A'=A_n$ by
Proposition \ref{flop} (3) again, a contradiction.
Thus $\beta'_i$ ($1\leq i\leq 6$) are
the only $g$-exceptional curves.
Now we investigate non-trivial fibers of $f'$.
Let $\gamma$ be a non-trivial fiber of $f'$. Then it holds that
\begin{equation}
\label{eq:March}
\text{$-K_{A'}\cdot \gamma=1$ and $L\cdot \gamma=0$.}
\end{equation}
If $\gamma$ is disjoint from all the flopped curves on $A'$,
then it holds also that 
$-K_{A}\cdot \gamma=1$ and $L\cdot \gamma=0$ on $A$
since $A\dashrightarrow A'$ is isomorphic near $\gamma$.
The equalities $-K_{A}=2H-E$ and $L=3H-2E$ show that
$H\cdot \gamma=2$ and $E\cdot \gamma=3$.
This means that the image of $\gamma$ on $B'$ 
is an irreducible tri-secant conic.
If $\gamma$ intersect some flopped curve on $A'$,
then $\gamma$ intersect only one flopped curve $\beta'$
at one point by Proposition \ref{prop:Cd1} (3).
Then, by applying Proposition \ref{prop:Atiyah} to the flop 
$A\dashrightarrow A'$,
the equalities (\ref{eq:March}) and $L\cdot \beta'=-1$ imply that
$-K_{A}\cdot \gamma=1$ and $L\cdot \gamma=1$ on $A$.
Thus we have  
$H\cdot \gamma=1$ and $E\cdot \gamma=1$.
Since $\gamma$ intersect the flopping curve on $A$
corresponding to $\beta'$,
the image of $\gamma$ on $B$ 
is a line as desired.\\
{\bf Step 4.}
We finish the proof 
by disproving the case where $f'$ contracts a divisor $E'$ to a point. 
 
By \cite{ray}, $f'$ is the blow-up at a point $b$ of $B'$ and
satisfies one of the following $E_2$--$E_5$:\\
$E_2:$ $b$ is a smooth point of $B'$. $E'\simeq \mP^2$ and 
$-K_{A'|E'}=\sO_{\mP^2}(2)$.\\
$E_3:$ $B'$ is analytically isomorphic to
$\{xy+zw=0\}\subset \mC^4$ near $b$. $E'\simeq \mP^1\times \mP^1$ and
$-K_{A'}={f'}^*(-K_{B'}¥)-E'$.\\ 
$E_4:$ $B'$ is analytically isomorphic to
$\{xy+z^2+w^3=0\}\subset \mC^4$ near $b$. 
$E'$ is a singular quadric surface and $-K_{A'}={f'}^*(-K_
{B'}¥)-E'$.\\ 
$E_5:$ $b$ is a $\frac 12(1,1,1)$-singularity. $E'\simeq \mP^2$ and
$-K_{A'}={f'}^*(-K_{B'})-\frac 12 E'$.

For the strict transform $q$ of a general tri-secant conic of $C$,
it holds that $-K_{A'}\cdot q=1$.
Therefore the case $E_2$ does not occur.
If $f'$ is of type $E_3$ or $E_4$,
then, by a similar consideration to Step 3,
we see that $B'$ is a (singular) quintic del Pezzo $3$-fold.
On the other hand, by $-K_{A'}={f'}^*(-K_{B'})-E'$ and ${E'}^3=2$,
we have $(-K_{B'})^3=(-K_{A'})^3+2=(-K_{A})^3+2=16$,
a contradiction.
If $f'$ is of type $E_5$, then,
by $-K_{A'}={f'}^*(-K_{B'})-\frac 12 E'$ and ${E'}^3=4$,
we have $(-K_{B'})^3=(-K_{A'})^3+\frac 12=(-K_{A})^3+\frac 12=\frac{29}{2}$.
By the classification of $\mQ$-Fano $3$-folds 
with only $\frac 12(1,1,1)$-singularities
(see \cite{Sa95}, \cite{Sa96}) and $(-K_{B'})^3=\frac{29}{2}$,
the possible Fano index of $B'$ is $\frac 12$, namely,
$2(-K_{B'})=\widehat{L}$.
If the Fano index of $B'$ is $\frac 12$, then
it holds that $2(-K_{A'})+{E'}=2{f'}^*(-K_{B'})=L=3H-2E$.
By $-K_{A'}=2H-E$, we have $E'\sim H$, a contradiction.
\end{proof}

\begin{rem}
It is possible to prove 
the existence of $\sS$ as in the statement of
Proposition \ref{prop:tworay}
but we do not prove this since
we do not need this in the sequel.
We only mention that $C_{\varphi}\in \sS$,
where $C_{\varphi}$ is the unique closed orbit of $G$-action on $B$. 
\end{rem}



\subsection{The correspondence between lines on $A$ and lines on $A'$}
\label{subsection:trigonality}~

The contents of this subsection is
presented also in \cite[the proof of Lemma 4.0.5]{TZb};
here we need a very detailed version of it for later usage.
Let $C$ be a general sextic normal rational curve on $B$
and we consider the diagram (\ref{eq:(A)}) as in Proposition \ref{prop:tworay}.
Denote by $\beta'_i$ the strict transform on $A$ of $\beta_i$.
Since $f'\colon A'\to B$ is also the blow-up of $B$
along a general sextic normal rational curve,
we can define the notion of lines on $A'$.
For simplicity of notation,
we denote by the same notation 
the strict transforms on $A'$ of curves and divisors on $A$.

\begin{prop}
\label{prop:AA'}
There exists a natural one to one correspondence between 
lines on $A$ and on $A'$.
\end{prop}

\begin{proof}
Let $l$ be a line on $A$. 
Since $E'=4H-3E$,
we have $E'\cdot l=1$ on $A$.
If $l$ is disjoint from any $\beta'_i$,
then we have $-K_{A'}\cdot l=1$ and $E'\cdot l=1$ on $A'$, thus
$l$ is a line on $A'$.
Assume that $l$ intersects some flopping curve of $A\dashrightarrow A'$.
By the classification of lines on $A$ (Proposition \ref{lineeA}),
there are two cases:
\begin{enumerate}[(a)]
\item
$l$ is the strict transform of a line on $B$
intersecting both $C$ and one bi-secant line $\beta_i$
outside $C\cap \beta_i$.

By Corollary \ref{extraline},
$l$ is the strict transform of
$\alpha_{i1}$ or
$\alpha_{i2}$.
Since $E'\cdot \beta'_i=-2$ on $A$,
it holds that $-K_{A'}\cdot l=1$ and $E'\cdot l=-1$ on $A'$
by Proposition \ref{prop:Atiyah}.
As in the proof of Proposition \ref{prop:tworay} (2),
$l$ is a fiber of $E'\to C'$.
Hence the union of $l\cup \beta''_i$, where
$\beta''_i$ is the flopped curve corresponding to $\beta'_i$,
is a line on $A'$ of type (ii) as in Proposition \ref{lineeA}.

In this case, $l$ corresponds to the line $l\cup \beta''_i$ on $A'$.
\item
$l$ is the union of the strict transform $\beta'_i$ of
one $\beta_i$ and a fiber $\zeta_{ij}$ of $E$ 
over one point $p_{ij}$ of $C\cap \beta_i$.

Note that $-K_A\cdot \zeta_{ij}=1$ and $E'\cdot \zeta_{ij}=3$ on $A$.
By $E'\cdot \beta'_i=-2$ on $A$ 
and Proposition \ref{prop:Atiyah},
it holds that $-K_{A'}\cdot \zeta_{ij}=1$ and $E'\cdot \zeta_{ij}=1$ on $A'$.
Therefore $\zeta_{ij}$ is a line on $A'$.
Moreover $f'(\zeta_{ij})$ is a line on $B$
intersecting $C'$ and the bi-secant line $\widehat{\beta}_i$
outside $C'\cap \widehat{\beta}_i$,
where $\widehat{\beta}_i$ is the image of $\beta''_i$ by $f'$.

In this case, $l$ corresponds to the line $\zeta_{ij}$ on $A'$.
\end{enumerate}
Thus, in any case, a line on $A$ corresponds to the unique line
on $A'$ and vice versa by symmetry of the diagram (\ref{eq:(A)}).
\end{proof}
We denote by $\sH'_1$ the curve obtained from $C'$ as a triple cover
as in Proposition \ref{primaC}.
The curve $\sH'_1$ has an ineffective theta characteristic $\theta'$
as in Proposition \ref{prop:sopratheta}.
\begin{prop}
\label{prop:HH'}
$(\sH_1,\theta)$ and $(\sH'_1,\theta')$ are isomorphic to each other
as spin curves.
\end{prop}

\begin{proof}
Since $\sH_1$ and $\sH'_1$ are the Hilbert schemes
of lines on $A$ and $A'$ respectively,
we can naturally identify $\sH_1$ and $\sH'_1$
by Proposition \ref{prop:AA'}.
Moreover, we can identify also $\theta$ and $\theta'$
since the strict transforms of 
two general intersecting lines on $A$
also intersect on $A'$ and vice versa, and 
the theta characteristics
are defined by the intersection of lines (see \ref{subsubsection:theta}).
\end{proof}

As we reviewed in the introduction,
the natural rational map
$\pi_{\sS^{+}_{4}}\colon  
\sH\dashrightarrow \sS^{+}_{4}$,
$C\mapsto (\sH_1,\theta)$
is the composite of the rational maps
$p_{\sS^{+}_{4}}\colon \sH\dashrightarrow
\widetilde{\sS}^{+}_{4}$
and $q_{\sS^{+}_{4}}\colon \widetilde{\sS}^{+}_{4}\dashrightarrow
{\sS}^{+}_{4}$,
where a general fiber of $p_{\sS^{+}_{4}}$ is 
a $\PGL_2$-orbit in $\sH$ and $q_{\sS^{+}_{4}}$ is 
birational or of degree two.
From Proposition \ref{prop:HH'}, we immediately obtain the following:
\begin{cor}
\label{cor:HH'}
The rational map $q_{\sS^{+}_{4}}$ exchanges
the classes of $C$ and $C'$ on 
$\widetilde{\sS}^{+}_{4}$.
\end{cor}

\section{The rationality proof}
\label{section: unirazionaledue}

\subsection{$\sH$ is birational to $(\mP^2)^6/\mathfrak{S}_6$.}~

By Proposition \ref{prop:irred}, 
$\sH$ is an irreducible $12$-dimensional variety.
The $G$-action on $B$ induces the $G$-action on $\sH$.
We construct a $G$-equivariant birational morphism
$\Theta\colon \widetilde{U}_0\to (\mP^2)^6/\mathfrak{S}_6$,
where
$\widetilde{U}_0$ is the open subset of $\sH$ consisting of 
(possibly reducible)
sextic curves with exactly six different bi-secant lines.

We remind the readers that the Hilbert scheme $\sH^B_1$ of lines on $B$
is isomorphic to $\mP^2$ and the Hilbert scheme $\sH_1$ of lines on $A$
is contained in the universal family $\mP$ of lines on $B$.
We know that the restriction to 
$\sH_{1}\subset\mP$ of the natural morphism 
$\pi\colon\mP\rightarrow\sH^{B}_{1}\simeq \mP^{2}$  is the morphism $\pi_{\mid \sH_{1}}\colon
\sH_1\to M$, where, by Corollary \ref{orabisecanti},
$M=\pi(\sH_1)$ is a plane nodal sextic whose nodes are the 
points $\beta_{1}, \ldots,
\beta_{6}\in\mP^{2}$ corresponding to the six bi-secant lines
$\beta_{1},\ldots,\beta_{6}\subset B$ of $C$. 

Then it remains defined a $G$-equivariant morphism 
$$
\Theta\colon\widetilde{U}_0 \rightarrow (\mP^2)^6/\mathfrak{S}_6,\,\,
C_{}\mapsto (\beta_{1},\ldots,\beta_{6}).
$$

\begin{lem}\label{pointsonlydominant} The morphism 
    $\Theta$ is dominant.
\end{lem}

\begin{proof} 
We need to prove that,
for six general lines on $B$,
there is a sextic rational curve on $B$ having them as its bi-secant
lines.
Let $\sG\subset \widetilde{U}_0$ be the divisor whose general
point parameterizes 
the union of a smooth quintic rational curve $C_5$ and a line $l$
such that they intersect simply at only one point.
First we prove that the restriction of $\Theta$ on
$\sG$ is dominant over the divisor $\sG'$ 
consisting of $6$-ples with three collinear points.
To show this,
let $l, l_1,l_2,l_3$ be four general lines on $B$ and
$m_1, m_2, m_3$ three general lines intersecting $l$.
We have only to prove there exists a $C_5$ such that
$C_5\cap l\not=\emptyset$, 
$l_1,l_2,l_3$ are three bi-secants of $C_5$ and
$m_1, m_2, m_3$ intersect both $C_5$ and $l$.
Consider the projection of $B$ from $l_1$.
Recall that the divisor $T_{l_1}$ swept by lines intersecting $l_1$
is mapped to a twisted cubic $\gamma$ on $Q$ by
Proposition \ref{projline}. The lines
$l, l_2,l_3, m_1,m_2,m_3$ are mapped to lines 
$l', l'_2,l'_3, m'_1,m'_2,m'_3$
intersecting $\gamma$.

Let $S$ be the smooth hyperplane section of $Q$
spanned by $l'_2$ and $l'_3$.
Note that $S$ is $\mP^1\times \mP^1$, and
$l'_2$ and $l'_3$ belong to the same ruling.
Let $n$ be a line in the other ruling.
Then by a simple dimension count,
there exists a twisted cubic $C'\sim 2n+l'_2$
passing through $5$ points 
$l'\cap S, m'_1\cap S,m'_2\cap S,m'_3\cap S$
and one point in $\gamma\cap S$.
The strict transform on $B$ of $C'$ is
a $C_5$ such that $C_5\cup l$ is mapped by $\Theta$ to
$(l_1,l_2,l_3,m_1, m_2, m_3)$.

We have proved that the divisor $\sG'$ in $(\mP^2)^6/\mathfrak{S}_6$
given by six lines such that three of them
intersect a line is dominated by $\sG$.
This is sufficient for the dominancy of $\Theta$.
Indeed,
for a general $C\in \widetilde{U}_0$, the $6$ points $\beta_1,\dots,\beta_6$
are in a general position
by Corollary \ref{cor:gen}, 
hence
$\Ima \Theta$ is not contained in $\sG'$.
Therefore, by the irreducibility of 
$(\mP^2)^6/\mathfrak{S}_6$ the claim follows.
\end{proof}

\begin{thm}\label{pointsonly} The morphism 
    $\Theta$ is birational.
\end{thm}

\begin{proof} 
Since $\dim \sH=\dim (\mP^{2})^6/\mathfrak{S}_6=12$,
it suffices to show that $\Theta$ is generically injective.

Let $\sH^o$ be the open set of $\sH$
consisting of sextic normal rational curves $C$
which satisfy all the following conditions: 
\begin{enumerate}[(a)]
\item
$C$ has exactly six different bi-secant lines $\beta_1,\dots,\beta_6$
(Corollary \ref{orabisecanti}).
\item
$\beta_1,\dots,\beta_6\in \mP^2$ are in a general position
(Corollary \ref{cor:gen}).
\item
For the strict transform $\beta'_i$ on $A$,
it holds that
$\sN_{\beta'_i/A}=\sO_{\mP^1}(-1)^{\oplus 2}$
(Proposition \ref{-1-1}).
\item
There are two lines $\alpha_{i1}$ and $\alpha_{i2}$ intersecting 
both $C$ and $\beta_i$ outside $C\cap \beta_i$
(Corollary \ref{extraline}).
\item
For a general line
$\alpha$ intersecting $\beta_i$,
there exist four lines $\gamma_1,\dots,\gamma_4$ 
different from $\beta_i$ and
intersecting both $C$ and $\alpha$ 
(Proposition \ref{prop:Cd} (5)).

To state the condition (f), we have the following remark:
note that it is possible to define the diagram
(\ref{eq:(A)}) as in Proposition \ref{prop:tworay}
since $C$ is normal and $C$ has only a finite number of bi-secant lines by (a).
Let $C'$ be as in Proposition \ref{prop:tworay}.
If we consider the blow-up $f'\colon A'\to B$ as the starting
point of the diagram $(\ref{eq:(A)})$,
then we obtain the diagram ending with $f\colon A\to B$. 
By this symmetry,
if $C$ is a general sextic normal rational curve,
then so is $C'$. 

\item
$C'$ is also contained in $\sH^o$.

\end{enumerate}

Let $\beta_{1}\ldots, \beta_{6}$
be general six lines on $B$.
By Lemma \ref{pointsonlydominant}, we may assume that
there is 
a sextic normal rational curve $C$
such that $C\in \sH^{o}$ and 
$\beta_{1}\ldots, \beta_{6}$
are bi-secant lines of $C$.
Let $\alpha_{ij}$ be as in the property (d).
\begin{cla}
\label{cla:main}
$\alpha_{ij}$ does not depend on $C$,
namely, if 
$\Gamma$ is another sextic normal rational curve
such that $\Gamma\in \sH^{o}$ and 
$\beta_{1},\ldots, \beta_{6}$
are also bi-secant lines of $\Gamma$,
then
$\alpha_{i1}$ and $\alpha_{i2}$ intersect 
$\Gamma$ outside $\Gamma\cap \beta_i$.
\end{cla}

\begin{proof}[Proof of the claim]
We take a general line
$\alpha$ intersecting $\beta_i$.
By the property (e),
there exist four lines $\gamma_1,\dots,\gamma_4$ 
different from $\beta_i$ and
intersecting both $C$ and $\alpha$. 
Let $\alpha'$ and $\beta'_i$ be
the strict transforms of $\alpha$ and $\beta_i$ on $A$.
We consider the six lines
$\alpha'_{i1}$, $\alpha'_{i2}$ and
$\gamma'_1,\dots,\gamma'_4$ on $A$
which are
the strict transforms on $A$ of 
the six lines 
$\alpha_{i1}$, $\alpha_{i2}$ and
$\gamma_1,\dots,\gamma_4$.

It holds that 
\begin{equation}
\label{eq:gamma}
\gamma'_1+\cdots+\gamma'_4=
(\pi_{|\sH_1})^*(M(\alpha)_{|M})-l_{i1}-l_{i2}
\end{equation}
and
\begin{equation}
\label{eq:alpha}
\alpha'_{i1}+\alpha'_{i2}\sim
(\pi_{|\sH_1})^*(M(\beta_i)_{|M})-(\delta-l_{i1})-(\delta-l_{i2}),
\end{equation}
where $l_{ij}$ are lines on $A$ as in Proposition \ref{lineeA} (ii).
Summing up these two equalities, 
we obtain
\begin{equation}
\label{eq:sum}
\alpha'_{i1}+\alpha'_{i2}+
\gamma'_1+\cdots+\gamma'_4
\sim 
(\pi_{|\sH_1})^*\sO_M(2)-2\delta=2\theta\sim K_{\sH_1},
\end{equation}
namely,
$\alpha'_{i1}+\alpha'_{i2}+
\gamma'_1+\cdots+\gamma'_4$ is a hyperplane section of $\sH_1\subset
\mP^4$.

Let $S\to \mP^2$ be the blow-up at six points 
$\beta_{1},\ldots, \beta_{6}$.
Let $\lambda\subset S$ be the total transform of a line on $\mP^2$
and $\varepsilon_i$ the exceptional curve over the point $\beta_i$.
Since $\varepsilon_{i|\sH_1}=l_{i1}+l_{i2}$,
the equality (\ref{eq:gamma}) implies that
\[
\gamma'_1+\cdots+\gamma'_4\sim (\lambda-\varepsilon_i)_{|\sH_1}.
\]
Thus, by the equality (\ref{eq:sum})
and
$K_{\sH_1}\sim
(3\lambda-\sum_{j=1}^{6} \varepsilon_j)_{|\sH_1},$ 
it holds that
\[
\alpha'_{i1}+\alpha'_{i2}\sim
\{(3\lambda-\sum_{j=1}^{6} \varepsilon_j)-(\lambda-\varepsilon_i)\}_{|\sH_1}= 
\{2\lambda-(\varepsilon_1+\cdots+\check{\varepsilon}_i+\cdots \varepsilon_6)\}_{|\sH_1}.
\]
Since $\sH_{1}$ is not hyperelliptic,
$\alpha'_{i1}+\alpha'_{i2}$ does not move,
thus
$\alpha'_{i1}+\alpha'_{i2}$ is cut out by
the strict transform of the unique conic $g_i$ on $\mP^2$
passing through $\beta_1,\dots,\check{\beta}_i,\dots,\beta_6$
(note the property (b)).
On the other hand, $\alpha_{i1}$, $\alpha_{i2}$ belong to
$M(\beta_i)$.
Thus $\alpha_{i1}$ and $\alpha_{i2}$ are exactly two points
of the intersection $g_i\cap M(\beta_i)$.
In particular this does not depend on $C$.
\end{proof}

Now we prove that $\Theta_{|\sH^o}$ is of degree one.
By contradiction 
assume that $\Gamma$ is a sextic rational curve
different from $C$ 
such that $\Gamma\in \sH^{o}$ and 
$\beta_{1}\ldots, \beta_{6}$
are bi-secant lines of $\Gamma$.
By the remark just before the property (f),
we can consider the diagram (\ref{eq:(A)})
as in Proposition \ref{prop:tworay} for $C$
and we use the notation there freely.
Let $\Gamma'$ be the strict transform of $\Gamma$
on $A$.
For simplicity of notation,
we denote by the same symbol the strict transforms on $A$ and $A'$ 
of curves on $B$.
On $B$ on the right hand side in the diagram $(\ref{eq:(A)})$,
let $\widehat{\Gamma}$ be the strict transform of $\Gamma$
and $\hat{\beta}_i$ 
the image of the flopped curve corresponding to $\beta_i$.

Since $\deg \Gamma=6$, we have $H\cdot \Gamma'=6$ on $A$.
By Proposition \ref{flop} (2),
it holds that $H\cdot \Gamma'\geq 6$ on $A'$.
Since $L$ is nef on $A'$ and $\Gamma'$ is not a fiber of 
$A'\to B$ by Proposition \ref{prop:tworay} (2),
it holds that $L\cdot\Gamma'\geq 1$.
Thus it holds $-K_{A'}\cdot \Gamma'\geq 4$ by 
$-2K_{A'}=H+L$ (cf. $(\ref{eq:(B)})$).
By Proposition \ref{flop} (1), it holds
$-K_{A}\cdot \Gamma'\geq 4$.
On the other hand,
$-K_B\cdot \Gamma=12$ on $B$ on the left hand side 
in the diagram $(\ref{eq:(A)})$.
Therefore,
since $-K_A=f^*(-K_B)-E$,
$\Gamma$ intersects $C$ at less than or equal to
$8$ points.
Thus, by the pigeon principle,
for at least two bi-secant lines of $C$, say,
$\beta_1$ and $\beta_2$,
$\Gamma$ passes through at most one 
of $p_{11}, p_{12}, t_{11}, t_{12}$ and
one of $p_{21}, p_{22}, t_{21}, t_{22}$,
where $\{p_{i1},p_{i2}\}:=C\cap \beta_i$ and
$t_{ij}:=C\cap \alpha_{ij}$
($i=1,2$, $j=1,2$).

This implies that 
$\hat{\beta}_1$ and $\hat{\beta}_2$ are
at least $3$-secant lines of $\widehat{\Gamma}$.
Indeed,
$\alpha'_{ij}$ on $A'$ is a fiber of $f'$
intersecting $\hat{\beta}_i$ by Proposition \ref{prop:tworay} (2).
In particular, this implies that
$\hat{\beta}_i$ ($i=1,2$) intersects
$\widehat{\Gamma}$ at more than or equal to $3$ points counted with
multiplicities
(if $\hat{\beta}_i$ passes 
through a singular point of $\widehat{\Gamma}$,
then we regard 
$\hat{\beta}_i$ as a multi-secant line of  
$\widehat{\Gamma}$).

Now we show that $\deg \widehat{\Gamma}\leq 6$.
Indeed, define the non-negative integer $a$
by the equation $-K_A\cdot \Gamma'=12-a$,
equivalently, $C$ intersects $\Gamma$ on $B$ on the left hand side
at $a$ points 
counted with multiplicities.
Then $\Gamma'$ intersects $\beta'_1,\dots,\beta'_6$
at more than or equal to $12-a$ points
depending on the common intersection points of $C$, $\Gamma$ and $\beta_i$.
This implies that $H\cdot \Gamma'\geq 6+12-a$ on $A'$
by Proposition \ref{prop:Atiyah}.
By $(\ref{eq:(B)})$ in Proposition \ref{prop:tworay}, 
we have $L\cdot \Gamma'\leq 2(12-a)-(18-a)=6-a\leq 6$ on $A'$. Thus $\deg \widehat{\Gamma}\leq 6$.
 
Consider the projection $B\dashrightarrow Q$
from the line $\hat{\beta}_1$ (see Proposition \ref{projline}).
Then the degree of the image $\widehat{\Gamma}'$ 
of $\widehat{\Gamma}$ is at most $3$
since $\deg \widehat{\Gamma}\leq 6$ and $\hat{\beta}_1$ is at least a $3$-secant line of $\widehat{\Gamma}$.
The lines $\hat{\beta}_i$ $(i=1,\dots,6)$ are the bi-secant lines of $C'$.
It holds that $\hat{\beta}_1\cap \hat{\beta}_2=\emptyset$
since $C'\in \sH^o$ by the property (f).
Thus the image of $\hat{\beta}_2$ on $Q$
is at least $3$-secant lines of $\widehat{\Gamma}'$.
If $\deg \widehat{\Gamma}'=1,2$, then this is impossible.
If $\deg \widehat{\Gamma}'=3$, then $\widehat{\Gamma}'$ is a twisted cubic 
curve since a plane cubic curve does not exist on $Q$.
Thus, again, $\widehat{\Gamma}'$ cannot have a $3$-secant line. 
\end{proof}

\subsection{Birational model of $\sS^+_4$}
\label{subsection:bir}~

Recall that
$\widetilde{U}_0$ is the open subset of $\sH$ consisting of 
sextic curves with exactly six different bi-secant lines.
Let $\widetilde{U}_1\subset \widetilde{U}_0$ be the open subset 
such that $\Theta$ is an isomorphism on $\widetilde{U}_1$.
Clearly $\widetilde{U}_1$ is $G$-invariant.
Let $U_1$ be the image of $\widetilde{U}_1$ on 
$(\mP^2)^6/\mathfrak{S}_6$.
Let $\widehat{U}\subset (\mP^2)^6$ be the set of stable 
ordered six points with respect to
the symmetric linearization of the action of $\PGL_3$, 
more explicitly,
the set of ordered six points such that
no two points coincide, or
no four points are collinear (see \cite[p.23, Theorem 1]{DO}).
By this explicit description,
we see that $\widehat{U}$ is $\mathfrak{S}_6$-invariant.
Note that the geometric quotient $\widehat{U}/G$ exists. Indeed,
let $\sL$ be the restriction of 
the $\PGL_3$-linearized line bundle to $\widehat{U}$. 
By restricting the $\PGL_3$-action to the $G$-action, 
$\sL$ is also $G$-linearized. 
We claim that
$\widehat{U}$ is the set of $G$-stable points.
Indeed, let $x\in \widehat{U}$ be a point.
The stabilizer group of $x$ for the $G$-action is finite (actually trivial)
since so is for the $\PGL_3$-action.
There exists a $\PGL_3$-invariant
section $s$ of some multiple of $\sL$ such that
$s(x)\not =0$ and $\PGL_3 \cdot x$ is closed in 
$\widehat{U}_s:=\{y\in \widehat{U}\mid s(y)\not =0\}$.
Since $G\subset \PGL_3$ is a closed subgroup,
the same is true for $G$.
 
Set ${U}_2=\widehat{U}/\mathfrak{S}_6\subset (\mP^2)^6/\mathfrak{S}_6$.
Since the $G$-action and $\mathfrak{S}_6$-action commutes,
$U_2/G$ also exists and $U_2/G\simeq (\widehat{U}/G)/\mathfrak{S}_6$.
Let $U'_3$ be the open subset of $U_1\cap U_2$
such that 
$C\in \Theta^{-1}(U'_3)$ is a sextic normal rational curve.
Note that, if $C\in \Theta^{-1}(U_1)$ is a sextic normal rational curve,
then we can define
the diagram (\ref{eq:(A)}) as in Proposition \ref{prop:tworay}
for $C$ since $C\in \Theta^{-1} (U_1)$ has only a finite number of bi-secant lines. Let $C'$ be as in Proposition \ref{prop:tworay}.
If $C$ is a general sextic normal rational curve, then
so is $C'$ by the symmetry of the diagram (\ref{eq:(A)}).
Thus $\Theta(C)\in U'_3$ with $C\in \Theta^{-1}(U'_3)$ such that
$C'\not \in \Theta^{-1} (U'_3)$
form a proper closed subset of $U'_3$, which we denote by $T$.
Set $U_3:=U'_3-T$, namely, $U_3$ is the biggest open subset of $U_1\cap U_2$
such that
$C\in \Theta^{-1} (U_3)$ is a sextic normal rational curve,
and
the center $C'$ of $f'\colon A'\to B$ also belongs to $\Theta^{-1} (U_3)$.
It is easy to see that $U_3$ is $G$-invariant
since the diagram  (\ref{eq:(A)}) is $G$-equivariant. 
Then by Corollary \ref{cor:HH'} and 
Theorem \ref{pointsonly}, the involution associated to the
map $q_{\sS^{+}_{4}}\colon \widetilde{\sS}^{+}_{4}\dashrightarrow
{\sS}^{+}_{4}$ is translated to an involution $J$ on $U_3/G$ satisfying
$J\colon \Theta(C) \mapsto \Theta(C')$ 
since $\widetilde{\sS}^{+}_{4}$
birationally parameterizes $G$-orbits in $\sH$.

We can sum up the above discussion into the following:
\begin{prop}
\label{prop:U_3}
$\sS^+_4$ is birational to
$(U_3/G)/J$.
\end{prop}
We investigate the variety $(U_3/G)/J$
relating it with 
the following classically well-studied variety:
\[
\sY:= (\mathbb P^{2})^{6} /\!/\PGL_3,
\] 
where the GIT-quotient is taken with respect to the symmetric linearization of the action of $\PGL_3$ (\cite[p.7, Proposition 1]{DO}). 
This is a compactification of the moduli space of ordered six points on 
$\mP^2$.
Note that
there exists a natural morphism
$U_3/G\to \sY/\mathfrak{S}_6$
since $G$-action on $(\mP^2)^6$ commutes with
$\mathfrak{S}_6$-action on $(\mP^2)^6$.

\subsection{A lifting of the association map on $(\mP^2)^6/\mathfrak{S}_6$ 
modulo $G$}~

We show that $J$ is a lifting of the classical association map on
$\sY/\mathfrak{S}_6$.

By \cite[p.37, Example 4]{DO},
there exists an involution $j'$ on $\sY$ called the 
{\em (ordered) association map}.
We do not give the definition of $j'$ but only describe it
on the open subset of $\sY$ which parameterizes
ordered six points in general positions (see \cite[p.118--120]{DO}). 

Let $\Sigma\subset \mP^3$ be a smooth cubic surface
and $\sigma\colon \Sigma\to \mP^2$
be the blow-up of $\mP^2$ at six points $p_1,\dots,p_6$.
We consider ordered sets of six lines on $\Sigma$,
equivalently, ordered sets of six points on $\mP^2$,
while till now, we have considered only unordered sets of six points
on $\mP^2$.
The 27 lines on $\Sigma$ can be grouped
into three ordered subsets:
$$(l_1,...,l_6),\,\, (l'_1,...,l'_6),\,\, (m_{ij}) \,
 (1\leq i<j\leq 6),$$ 
where the lines $l_i$ are the exceptional lines $\sigma^{-1}
(p_i)$, the lines $l'_i$ are the strict transforms of
the conics $q_i\subset \mP^2$ passing through 
$p_1,\dots,\check{p}_i,\dots,p_6$,
and the lines
$m_{ij}$ are the strict transforms of the lines $\langle p_i, p_j\rangle$
joining the points $p_i$ and $p_j$.
The first two groups of lines $(l_1,...,l_6)$ and $(l'_1,...,l'_6)$
form a {\it double sixer}, which means that
$$l_j\cap l_j = \emptyset, \quad
\l'_i \cap \l'_j = \emptyset,\quad l_i\cap l'_j \neq \emptyset
\quad{\text{if and only if}}
\quad i\neq j.$$
Every set of 6 disjoint lines on $\Sigma$ can be included in
 a unique double sixer, from which
$\Sigma$ can be reconstructed uniquely.
 There are 36 double sixers of $\Sigma$. Every double sixer defines
two regular birational maps $\sigma: \Sigma \rightarrow \mP^2$, 
$\sigma': \Sigma
\rightarrow \mP^2$, each of which blows down one of the two sixes (sixtuples
of disjoint lines) of the double sixers. 
The association map $j'$ interchanges  
the two collections of ordered 6 points in $\mP^2$ given by $(\sigma(l_1),\dots,
\sigma(l_6))$
and  $(\sigma'(l'_1),\dots,
\sigma'(l'_6))$, namely,
it holds that
\begin{equation}
\label{eq:j}
j'\colon (\sigma(l_1),\dots,\sigma(l_6))\mapsto
(\sigma'(l'_1),\dots,\sigma'(l'_6)).
\end{equation}
We also remark that $j'$ fixes ordered six points on a conic.

Since the symmetric group $\mathfrak{S}_6$ acts   
on the quotient $\sY$ and
its action commutes with $j'$,
the map $j'$ descends to 
an involution $j$ on $\sY/\mathfrak{S}_6$.
The map $j$ is called the {\em (unordered) association map}.
 
\begin{prop}\label{nori}
The involution $J$ is a lifting of $j$.
\end{prop}

\begin{proof}
It suffices to check the assertion
at a general point of $U_3/G$.
Let $C\in \sH$ be a general point such that $\Theta(C)\in U_3$.
By definition of $\Theta$, it holds that
$\Theta(C)=(\beta_1,\dots, \beta_6)$,
where $\beta_1,\dots,\beta_6$ are six bi-secant lines of $C$.
Now we compute $\Theta(C')$.
By Corollary \ref{extraline},
there exist two lines $\alpha_{i1}$ and $\alpha_{i2}$ intersecting
a bi-secant line $\beta_i$ and $C$ outside $C\cap \beta_i$
($1\leq i \leq 6$).
Let 
$\alpha''_{i1}$ and $\alpha''_{i2}$ be the strict transforms
of $\alpha_{i1}$ and $\alpha_{i2}$
on $A'$.
Then $\alpha''_{i1}$ and $\alpha''_{i2}$ 
are the fibers of $E'$ through $E'\cap \beta''_i$
by Proposition \ref{prop:tworay} (2),
where 
$E'$ is the $f'$-exceptional divisor
and $\beta''_i$ is the flopped curve corresponding to
$\beta_i$.
Thus, by Corollary \ref{cor:node}, the two lines 
$\alpha''_{i1}\cup \beta''_i$ and $\alpha''_{i2}\cup \beta''_i$
on $A'$
correspond to the node $\widehat{\beta}_i$ of $M(C')$,
where $\widehat{\beta}_i$ is the image of $\beta''_i$ by $f'$.
Let 
$\beta'_i$, $\alpha'_{i1}$ and $\alpha'_{i2}$ be the strict transforms
of $\beta_i$, $\alpha_{i1}$ and $\alpha_{i2}$
on $A$.
By the proof of Proposition \ref{prop:AA'},
the two lines 
$\alpha''_{i1}\cup \beta''_i$ and $\alpha''_{i2}\cup \beta''_i$
on $A'$ correspond to
the two lines $\alpha'_{i1}$ and $\alpha'_{i2}$ on $A$.
By the proof of Claim \ref{cla:main},
$\alpha'_{i1}+\alpha'_{i2}$ on $\sH_1$ 
is the divisor cut 
by the strict transform of the unique conic passing 
through all the nodes of $M(C)$ except $\beta_i$.
Thus, by the above description of $j'$,
$J$ is a lifting of the association map $j$ on 
$\sY/\mathfrak{S}_6$. 
\end{proof}

\subsection{Rationality of the moduli space of double sixers on $\mP^2$}~

By Proposition \ref{prop:U_3},
we have only to show that $(U_3/G)/J$ is rational
for the proof of the main theorem.
Professor Igor Dolgachev kindly told us 
that a similar statement for
$\PGL_3$ is true, more precisely, 

\begin{prop}[A.~Coble]
\label{ortland} The quotient variety 
    $(\sY/\mathfrak{S}_6)/j$ is a rational variety.
\end{prop}

This is a classical result due to 
A.~Coble, which easily follows from \cite[p.19 and 37]{DO}.

This result is a bit subtle;
it is not known if the moduli space
$\sY/\mathfrak{S}_6$ of unordered six points on $\mP^2$
is rational or not.  

By the proof of
Proposition \ref{ortland},
we see that the degree of the map 
$\sY/\mathfrak{S}_6\to (\sY/\mathfrak{S}_6)/j$ 
is two, namely, we have
\begin{cor}
\label{cor:ortland}
The map $j$ is a non-trivial involution on $\sY/\mathfrak{S}_6$
whence so is $J$.
\end{cor}

\subsection{Proof of the rationality of $\sS^+_4$}~

The following diagram summerizes our construction above:
\begin{equation}
\label{eq:basicthree}
\xymatrix{\widetilde{U}_0\ar@{-->}[r]^{p_{\sS^+_4}}\ar@{-->}[d]_{\Theta}&
\widetilde{\sS}^+_4\ar@{-->}[r]^{q_{\sS^+_4}}\ar@{-->}[d]_{\mathrm{bir.}}&
{\sS}^+_4\ar@{-->}[d]_{\mathrm{bir.}}\\
U_3 \ar[r]& 
U_3/G \ar[d]_{\varrho}\ar[r]& 
(U_3/G)/J \ar[d]\\
& \sY/\mathfrak{S}_6  \ar[r]&
(\sY/\mathfrak{S}_6)/j} 
\end{equation}

After proving some lemmas we show that $(U_3/G)/J$ is a rational variety.

We consider the following diagram:
\begin{equation}
\label{eq:PGL_3}
\xymatrix{
\widehat{U} \ar[d]_{\widehat{\pi}_{\PGL_{3}}}\ar[r]& 
\widehat{U}/\mathfrak{S}_6 \ar[d]_{\pi_{\PGL_{3}}}\\
\widehat{U}/\PGL_3 \ar[r]^{h}\ar@{_{(}->}[d]& 
(\widehat{U}/\PGL_3)/\mathfrak{S}_6\ar@{_{(}->}[d]\\
\sY\ar[r] & \sY/\mathfrak{S}_6,}
\end{equation}
where 
recall that $\widehat{U}\subset (\mP^2)^6$ is the set of stable 
ordered six points with respect to
the symmetric linearization of the action of $\PGL_3$.

\begin{lem}\label{serre} 
The natural projection $\pi_{\PGL_{3}}$ is 
a principal fiber bundle of $\PGL_3$ over some non-empty 
open subset $W_1$ of $(\widehat{U}/\PGL_3)/\mathfrak{S}_6$.
\end{lem}
\begin{proof} 
By \cite[p.30, in the end of the proof of Theorem 2]{DO},
$\widehat{\pi}_{\PGL_{3}}$ is a principal fiber
bundle of $\PGL_3$.
We have seen that $\sY$ is
isomorphic to a quartic hypersurface in $\mP(1^5,2)$, hence
its degree $\sO_{\sY}(1)^4$ is equal to $2$.
By the equality (\ref{eqnarray:34}) in proof of Proposition \ref{ortland},
$\sY':=\sY/\mathfrak{S}_6$
is a hypersurface of degree $34$
in $\mP(2,3,4,5,6,17)$. 
Then
its degree $\sO_{\sY'}(1)^4$ is equal to $\frac{2}{6!}$.
Therefore the degree of the map $h$ in the diagram (\ref{eq:PGL_3})
is $6!$, which is equal to the order of $\mathfrak{S}_6$.
Hence $\mathfrak{S}_6$ acts trivially on
fibers of $\widehat{\pi}_{\PGL_{3}}$
over points in the open subset $W'_1$ of $\widehat{U}/\PGL_3$
where $h$
is \'etale.
By \cite[p,7, Proposition 0.2 and p.16, Proposition 0.9]{Mum2},
$\pi_{\PGL_{3}}$ is 
a principal fiber bundle of $\PGL_3$ over $W_1:=h(W'_1)$.
\end{proof}

Now we consider the $G$-action.
Let $\varrho\colon U_3/G\to \sY/\mathfrak{S}_6$
be the natural morphism.
Set $V_1:=\varrho^{-1} (W_1) \cap J (\varrho^{-1} (W_1))(\not =\emptyset)$.
By definition, $V_1$ is invariant under $J$.
Let $W_2:=\varrho(V_1)$ and $W'_2:=h^{-1} (W_2)$.

From the diagram (\ref{eq:PGL_3}) and the proof of
Lemma \ref{serre},
we obtain the following diagram:
\begin{equation}
\xymatrix{
{\widehat{\pi}_{\PGL_{3}}}^{-1}(W'_2)/G \ar[r]\ar[d]& 
{{\pi}_{\PGL_{3}}^{-1} (W_2)/G\supset V_1}\ar[d]_{\varrho'}\\
W'_2 \ar[r]\ar@{_{(}->}[d]& W_2\ar@{_{(}->}[d]\\
\sY\ar[r] & \sY/\mathfrak{S}_6}
\end{equation}

\begin{lem}\label{serrebis} The natural projection $\varrho'$
is a $\mP^5$-bundle.
\end{lem}
\begin{proof}  
    A fiber of $\varrho'$ is isomorphic to $\PGL_{3}/G$,
    which is isomorphic to $\mP^{5}$ by Proposition \ref{closure}.
\end{proof}    
Set $\sV:={\pi}_{\PGL_{3}}^{-1} (W_2)/G$.
We are going to find a sub $\mP^4$-bundle of $\rho'\colon \sV\to W_2$.

Recall that $\Omega\subset\sH^{1}_{B}=\mP^{2}$ is
the $G$-invariant conic (Proposition \ref{closure}),
and, for the symmetric bi-linear form $\widetilde{\Omega}$ 
associated to $\Omega$,
it holds that two lines $l$ and $m$ on $B$ intersect if and only if 
$\widetilde{\Omega}(l,m)=0$, where $l,m\in \sH^B_1$ are the points
corresponding to $l$ and $m$.
Let  $\sD'\subset (\mP^2)^6/\mathfrak{S}_6$
be the closure of the set of unordered six points 
two of which are polar with respect to $\widetilde{\Omega}$.
Clearly $\sD'$ is $G$-invariant.

\begin{lem}\label{intersecting}
    The locus $\sD'$ is an irreducible divisor of $(\mP^2)^6/\mathfrak{S}_6$.

For a general point $(l_1,\dots,l_6)\in \sD'$, it holds that
\begin{enumerate}[$(1)$]
\item
only two of six lines $l_1,\dots,l_6$ intersect on $B$, and
\item
six points $l_1,\dots,l_6\in \mP^2$ are in a general position.
\end{enumerate}  
\end{lem}

\begin{proof}

$\sD'$ is the image of the locus $\sD''$ defined by 
ordered six points $(l_1,\dots, l_6)\in(\mP^{2})^{6}$ such that
$\widetilde{\Omega}(l_5,l_6)=0$.
Since $(l_1,\dots, l_5)$ moves freely, the $5$-ples
$(l_1,\dots, l_5)$ are parameterized by $(\mP^2)^5$.
Once we fix $l_5$, the points $l_6$ are parameterized 
by the line $\widetilde{\Omega}(l_5,*)=0$. 
Then $\sD''$ is birational to a $\mP^1$-bundle over $(\mP^2)^5$.
In particular $\sD''$ is an irreducible divisor and so is $\sD'$.

Similarly, we can show that 
the sublocus in $\sD''$ consisting of $6$-ples
$(l_1,\dots, l_6)$ not satisfying (1) nor (2)  
is $4$-dimensional. Thus the latter assertion follows.  
\end{proof}

\begin{lem}\label{rationalnormal}
    A general point of $\sD'$ is the image
    by $\Theta$ of a sextic normal rational curve $C$
such that it is possible to define the diagram $(\ref{eq:(A)})$
as in Proposition $\ref{prop:tworay}$ for $C$.
\end{lem}

\begin{proof} 
Once we show that 
a general point of $\sD'$ is the image
    by $\Theta$ of a sextic normal rational curve $C$,
then it is possible to define the diagram $(\ref{eq:(A)})$
as in Proposition $\ref{prop:tworay}$ for $C$
since such a $C$ has only a finite number of bi-secant lines.

First we show that a general point of $\sD'$
is in the image of $\Theta$.
Let $H$ be a smooth hyperplane section of $B$. Then 
    $H$ is the blow-up of $\mP^2$ at four points. 
It is easy to show that, 
    if $C_o$ is the strict transform of a general smooth conic on $\mP^2$, 
then $C_o\in\sH$ by Proposition \ref{gradoduno6} (2),
$C_o$ has six b-secant lines,
and there exist three pairs of intersecting lines 
among the six bi-secant lines of $C_o$ . 
In particular $C_o\in \widetilde{U}_0$ and $\Theta(C_o)\in \sD'$.
Since $\Theta$ is dominant by Lemma \ref{pointsonlydominant}, 
a general point of $\sD'$
is also in the image of $\Theta$.


Take $C\in \widetilde{U}_0$ such that $\Theta(C)$ is a general
point of $\sD'$.
By contradiction assume that $C$ is contained in a hyperplane section
$H$ of $B$.
We may assume that $H$ is smooth since
even for the special sextic rational curve $C_o$ 
as above,
the hyperplane section containing it is smooth.  
Let $\beta_{1}$,\ldots , $\beta_{6}$
    be the six bi-secant lines of $C$.
By generality of $C$ and Lemma \ref{intersecting} (1),
we may assume that   
only $\beta_{5}$ and $\beta_{6}$ intersect. 
    Note that $\beta_i$ ($i=1,\dots,6$) are contained in $H$.
    We have a contraction $H\to \mP^2$
    of $\beta_1,\dots,\beta_4$ since they are disjoint,
    a contradiction since the image of $\beta_5$ and $\beta_6$ on $\mP^2$ are 
still $(-1)$-curves. 
\end{proof}

 Denote by $\sD\subset \sV={{\pi}_{\PGL_{3}}^{-1} (W_2)}/G$ 
 the image of
 $\sD'\cap {{\pi}_{\PGL_{3}}^{-1} (W_2)}$ by the quotient map.
\begin{lem}\label{invarioanceod d} $\sD\cap V_1\not =\emptyset$
and $\sD\cap V_1$ is invariant under $J$.
\end{lem}
\begin{proof} 
Recall the notation as in the subsection \ref{subsection:bir}.
By Lemma \ref{rationalnormal},
$\sD'$ intersects the image of $\Theta$.
Moreover, since $\sD'\subset (\mP^2)^6/\mathfrak{S}_6$
is a divisor, $\Theta$ is isomorphic over a general point of 
$\sD'$, namely, $U_1\cap \sD'\not =\emptyset$.
By Lemmas \ref{intersecting} (2) and \ref{rationalnormal},
$U'_3\cap \sD'\not =\emptyset$.
By Lemma \ref{rationalnormal}, 
we can define the diagram (\ref{eq:(A)}) as in 
Proposition \ref{prop:tworay} for $C\in \Theta^{-1}(U'_3\cap \sD')$.
    Let $\beta_1$ and $\beta_2$ be two intersecting bi-secant lines of $C$.
    Let $\beta'_i\subset A$ and $\beta''_i\subset A'$ ($i=1,2$) be the 
flopping and flopped curve corresponding to $\beta_i$ respectively.
Since $\beta'_1\cap \beta'_2\not =\emptyset$,
it holds that $\beta''_1\cap \beta''_2\not =\emptyset$ 
by Proposition \ref{prop:symm} (2).
Thus two of bi-secant lines of $C'$ which is the images of 
$\beta''_1$ and $\beta''_2$ intersect. 
This implies that
$\Theta(C')\in \sD'$.
By generality of $C$, we may assume that
$C'$ is also general.
In particular, we may assume that
$\Theta(C)\in \sD'\cap U_3$.
Then it holds that $J(\Theta(C))=\Theta(C')$ on $U_3/G$.
This implies that $\sD$ is invariant by the action of $J$.
\end{proof}

\begin{lem}\label{projectivebundle} The restriction to $\sD$ 
    of 
    $\varrho'\colon \sV \to W_2$
gives a $\mP^{4}$-bundle structure on $\sD$. 
\end{lem}
\begin{proof}
Let $m_1,\dots,m_6$ be six lines on $B$ such that
$(m_1,\dots,m_6)\in (\mP^2)^6/\mathfrak{S}_6$ is mapped to
a point $w$ of $W_2$. 
We show that the restriction of $\sD$ to
the fiber $F$ of $\varrho'$ over the point $w$
is isomorphic to $\mP^4$.
By Claim \ref{cla:doubly},
$G$ acts doubly transitively 
on the set of general unordered pairs of intersecting lines.
Thus, for any 
$1\leq i<j\leq 6$ and $1\leq k<l\leq 6$,
it holds that
\begin{eqnarray*}
\{(g_1(m_1),\dots,g_1(m_6))\in (\mP^2)^6/\mathfrak{S}_6\mid 
g_1\in \PGL_3, g_1(m_i)\cap g_1(m_j)\not =\emptyset\}
= \\
\{(g_2(m_1),\dots,g_2(m_6))\in (\mP^2)^6/\mathfrak{S}_6\mid 
g_2\in \PGL_3, g_2(m_k)\cap g_2(m_l)\not =\emptyset\}
\ \text{modulo}\ G
\end{eqnarray*}
since there exists $h\in G$
such that $h\{g_1(m_i),g_1(m_j)\}=\{g_2(m_k),g_2(m_l)\}$
by the double transitivity.
Therefore, a point of $F\cap \sD$ is the image of
a point $(g(m_1),\dots,g(m_6))\in (\mP^2)^6/\mathfrak{S}_6$, where 
$g\in \PGL_3$ and $\widetilde{\Omega}(g(m_5),g(m_6))=0$.
Now we choose a coordinate of $\mP^2$
such that $\Omega=\{x^2+y^2+z^2=0\}$.
Set $m_5=(a_1:a_2:a_3)$ and $m_6=(b_1:b_2:b_3)$.
Then 
$\widetilde{\Omega}(g(m_5),g(m_6))=0$
if and only if 
\begin{equation}
\label{eq:linear}
\begin{pmatrix}
a_1& a_2& a_3
\end{pmatrix}
{\empty^t\! g}g 
\begin{pmatrix}
b_1\\
b_2\\
b_3\\
\end{pmatrix}=0.
\end{equation}
Recall that by Proposition \ref{closure} the map 
$\PGL_3\to \mP_* H^0(\mP^2,\sO_{\mP^2}(2))\simeq \mP^5$
is defined by $g\mapsto {\empty^t \! g}g$,
where a conic on $\mP^2$ is identified with
a $3\times 3$ symmetric matrix.
 Since the condition (\ref{eq:linear})
is linear,
$F\cap \sD$ is a hyperplane in $F\simeq \mP^5$.
\end{proof}

\begin{lem}
\label{extend}
The involution $J$ on $V_1$ extends on $\sV$.
\end{lem}
 
\begin{proof}
By \cite[III Corollary 12.9]{Ha},
$\sF={\varrho'}_*\sO_{\sV}(\sD)$
is a locally free sheaf of
rank $6$ on
$W_2$ and 
${\varrho'}_*\sO_{\sV}(\sD)\otimes k(w)
\simeq H^0({\varrho'}^{-1}(w), \sO_{\sV}(\sD)_{|{\varrho'}^{-1}(w)})$
for $w\in W_2$. 
Consider the following diagram:
\begin{equation}
\label{eq:standard}
\xymatrix{{\varrho'}^*\sF \ar[d] \ar[r] & \sO_{\sV}(\sD)\ar[d]\\
{\varrho'}^*\sF_{|{\varrho'}^{-1}(w)}\ar[d] \ar[r] & \sO_{\sV}(\sD)_{|{\varrho'}^{-1}(w)}\ar[d]_{id}\\
H^0({\varrho'}^{-1}(w), \sO_{\sV}(\sD)_{|{\varrho'}^{-1}(w)})\otimes \sO_{{\varrho'}^{-1}(w)}
\ar[r] & {\sO_{\sV}(\sD)_{|{\varrho'}^{-1}(w)}}.}
\end{equation}
The map
$H^0({\varrho'}^{-1}(w), \sO_{\sV}(\sD)_{|{\varrho'}^{-1}(w)})
\otimes \sO_{{\varrho'}^{-1}(w)}
\to \sO_{\sV}(\sD)_{|{\varrho'}^{-1}(w)}$ is surjective since
${\varrho'}^{-1}(w)\simeq \mP^5$ and
$\sO_{\sV}(\sD)_{|{\varrho'}^{-1}(w)}\simeq \sO_{\mP^5}(1)$.
Thus, by the Nakayama lemma,
${\varrho'}^{*}\sF\rightarrow 
\sO_{\sV}(\sD)$ is surjective.
Then by \cite[II Proposition 7.12]{Ha} it remains
defined a morphism $\gamma\colon
\sV\rightarrow \mP(\sF)$ over 
$W_2$. Since $\gamma$ is
fiberwise an isomorphism then it is an isomorphism 
by the Zariski main theorem. 

Let $W^o$ be any open subset of $W_2$.
Since $\sD$ is invariant under the rational involution $J$,
it holds that 
$\Gamma(W^o, \sF)\simeq \Gamma(j(W^o),\sF)$,
which induces an isomorphism $\sF\simeq j^*\sF$.
Thus $J$ extends to the involution 
$\sV\simeq \mP(j^*\sF)\to \mP(\sF)=\sV$.
\end{proof}

We still denote by $J$ the extension of $J$ to $\sV$.
Set $\sR:=\sV/J$ and $\sW:=W_2/j$.
Now we can prove the main result:

\begin{thm}\label{minimalaim}
$\sR$ is a rational variety.
\end{thm}

\begin{proof}
The action of $J$ is trivial on the fiber of $\varrho'$
since $j$ acts non-trivially on 
$W_2$ by Corollary \ref{cor:ortland}.  
Thus $\varrho'$ descends to 
a $\mP^5$-bundle $p\colon \sR\to \sW$. 
Moreover, the sub-$\mP^4$-bundle $\sD$ of $\sV$
descends to a sub-$\mP^4$-bundle $\sT$ of $\sR$ since
it is invariant under $J$ by Lemma \ref{invarioanceod d}.  
Set $\sE:=p_*\sO_{\sR}(\sT)$. 
As in the proof of Lemma \ref{extend},
we can show that $\sR\simeq \mP(\sE)$.
In particular, $\sR$ is 
a locally trivial $\mP^5$-bundle over $\sW$.
Consequently $\sR$ is rational since so is $\sW$
by Proposition \ref{ortland}.
\end{proof}

\begin{cor} 
    ${\sS^{+}_{4}}$ is a rational variety.
\end{cor}
\begin{proof}
    It follows by Proposition \ref{prop:U_3} and Theorem 
    \ref{minimalaim} since $\sR$ is birationally equivalent to
    $(U_3/G)/J$.
\end{proof}



\end{document}